\documentclass{article}
\usepackage{maa-monthly}
\usepackage{comment}
\usepackage{lineno}

\usepackage{pstricks}
\psset{unit=3cm}
\usepackage{url}
\usepackage[section]{placeins}	
\usepackage{floatrow}
\usepackage{caption}
\usepackage{alltt}
\usepackage[greek,english]{babel}
\usepackage[utf8x]{inputenc}
\newrgbcolor{lightblue}{0.8 0.8 1}
\newrgbcolor{pink}{1 0.8 0.8}
\newrgbcolor{lightgreen}{0.8 1 0.8}
\newrgbcolor{lightyellow}{1 1 0.8}

\def\B{{\bf B\,}}   

\def\justification#1{\hfill{{\bf [#1]}\\}}

\def\qed{{\flushright{ \vspace{-15pt} Q.E.D.}\flushleft}}

\def\TarskiFiveSegmentFigure{%
\psset{unit=2.25cm}
\pspicture(2.5,0)(2.2, 0.95)
\pspolygon[fillstyle=solid,fillcolor=lightblue,linestyle=none]
(1,0.85)(2,0.15)(0.7,0.15)
\pspolygon[fillstyle=solid,fillcolor=lightblue,linestyle=none]
(3.5,0.85)(4.5,0.15)(3.2,0.15)
\qline(0.0,0.15)(2.0,0.15)
\qline(0.0,0.15)(1.0,0.85)
\psline[linestyle=dashed](1.0,0.85)(2.0,0.15)
\qline(1.0,0.85)(0.7,0.15)
\put(1,0.9){$d$}
\put(0.0,0){$a$}
\put(0.7,0){$b$}
\put(2.0,0){$c$}
\qline(2.5,0.15)(4.5,0.15)
\qline(2.5,0.15)(3.5,0.85)
\psline[linestyle=dashed](3.5,0.85)(4.5,0.15)
\qline(3.5,0.85)(3.2,0.15)
\put(3.5,0.9){$D$}
\put(2.5,0){$A$}
\put(3.2,0){$B$}
\put(4.5,0){$C$}
\psset{unit=3cm}
\endpspicture}

\def\SameSideFigure{%
\pspicture(0.55,0.45)(2,1.5)
\qline(0.7,1)(2,1)
\put(0.6,0.97){$p$}
\put(2.05,0.97){$q$}
\psdot(0.9,0.5)
\put(0.8,0.46){$a$}
\psdot(1.4,0.5)
\put(1.31,0.46){$b$}
\psline(1.8,1.5)(1.4,0.5)
\psline(1.8,1.5)(0.9,0.5)
\pscircle[fillstyle=solid,fillcolor=white](1.6,1){0.03}
\put(1.6,0.87){$y$}
\pscircle[fillstyle=solid,fillcolor=white](1.346,1){0.03}
\put(1.3,0.87){$x$}
\pscircle[fillstyle=solid,fillcolor=white](1.8,1.5){0.03}
\put(1.85,1.5){$c$}
\endpspicture
}

\def\OppositeSideFigure{%
\pspicture(0.55,0.45)(2,1.5)
\qline(0.7,1)(2,1)
\put(0.6,0.97){$p$}
\put(2.05,0.97){$q$}
\psdot(0.9,0.5)
\put(0.8,0.46){$a$}
\psdot(1.8,1.5)
\put(1.85,1.5){$b$}
\psline(1.8,1.5)(0.9,0.5)
\pscircle[fillstyle=solid,fillcolor=white](1.346,1){0.03}
\put(1.27,1.05){$x$}
\endpspicture
}

\def\CrossbarFigure{%
\pspicture(-0.5,0)(0.8,0.6)
\SpecialCoor
\qline(0,0)(0.45;60)
\qline(0,0)(0.7;120)
\qline(0,0)(0.9;30)
\psdot(0.5;30)
\psdot(0.75;30)
\psdot(0.5;120)
\psdot(0.3;120)
\psdot(0.292;60)
\psdot(0,0)
\qline(0.5;30)(0.3;120)
\put(-0.38,0.4){$J$}
\put(0.67,0.32){$K$}
\put(0.26,0.44){$X$}
\put(-0.3,0.2){$U$}
\put(-0.13,-0.05){$B$}
\put(0.47,0.17){$V$}
\put(0.09,0.3){$P$}
\put(-0.47,0.54){$A$}
\put(0.8,0.47){$C$}
\psline(0.5;120)(0.75;30)
\pscircle[fillstyle=solid,fillcolor=white](0.46;60){0.03}
\endpspicture
}

\def\InnerOuterPaschFigure{%
\pspicture(0,0.2)(3.6,1.6)
\pspolygon[fillstyle=solid,fillcolor=lightblue,linestyle=none]
(0,0.6)(1.1,1.4)(1.3,0.6)
\pspolygon[fillstyle=solid,fillcolor=lightblue,linestyle=none]
(2,0.3)(3.0,1.05)(3.5,0.3)
\psdot(0,0.6)
\put(0,0.47){$A$}
\pscircle[fillstyle=solid,fillcolor=white](1,0.6){0.03}
\put(1,0.47){$x$}
\qline(0,0.6)(0.97,0.6)
\psdot(1.3,0.6)
\put(1.35,0.55){$C$}
\qline(1.03,0.6)(1.3,0.6)
\psdot(1.1,1.4)
\put(1.1,1.45){$B$}
\qline(0,0.6)(1.1,1.4)  
\qline(1.1,1.4)(1.375,0.3) 
\psdot(1.375,0.3)
\put(1.43,0.25){$E$}
\psdot(0.52,0.98)
\put(0.43,1.03){$F$}
\qline(1.375,0.3)(1.015,0.584)  
\qline(0.52,0.98)(0.979,0.62)  
\psdot(2,0.3)
\put(1.95,0.16){$A$}
\psdot(3.5,0.3)
\put(3.45,0.16){$C$}
\qline(2,0.3)(2.47,0.3)
\qline(2.53,0.3)(3.5,0.3)
\pscircle[fillstyle=solid,fillcolor=white](2.5,0.3){0.03}
\put(2.45,0.19){$x$}
\psdot(2.7,1.5)
\put(2.65,1.57){$E$}
\qline(2.7,1.5)(3.5,0.3)  
\psdot(3.0,1.05)
\put(3.07,1.02) {$B$}
\qline(2.7,1.5)(2.5,0.33)  
\qline(2,0.3)(3.0,1.05)   
\psdot(2.566,0.722)
\put(2.4,0.7){$F$}
\endpspicture}

\def\ProclusBisectionFigure{%
\psset{unit=2cm}
\pspicture (2,0.9)(-1,-0.7) 
\pspolygon[fillstyle=solid,fillcolor=yellow,linestyle=none]
  (-0.25,-0.433)(0.5,0.866)(1.25,-0.433)(0.5,-0.176)
\psdot(0,0)
\psdot(1,0)
\psdot(0.5,0.866)
\psdot(-0.25,-0.433)
\psdot(1.25,-0.433)
\qline(-0.25,-0.433)(0.5,0.866)
\qline(0.5,0.866)(1.25,-0.433)
\qline(1.25,-0.433)(0,0)
\qline(-0.25,-0.433)(1,0)
\qline(0,0)(1,0)
\qline(0.5,0.866)(0.5,-0.176)
\pscircle[fillstyle=solid,fillcolor=white](0.5,-0.176){0.04}
\put(-0.2,-0.03){$A$}
\put(1.06,-0.03){$B$}
\put(-0.47,-0.48){$D$}
\put(1.3,-0.48){$E$}
\put(0.57,0.83){$C$}
\put(0.44,-0.38){$F$}
\psset{unit=3cm}
\endpspicture}
 
\def\ProclusBisectionFigureTwo{%
\psset{unit=2cm}
\pspicture (2,0.9)(-1,-0.7) 
\pspolygon[fillstyle=solid,fillcolor=yellow,linestyle=none]
  (0,0)(0.5,0)(0.5,0.866)(1.25,-0.433) 
\psdot(0,0)
\psdot(1,0)
\psdot(0.5,0.866)
\psdot(-0.25,-0.433)
\psdot(1.25,-0.433)
\qline(-0.25,-0.433)(0.5,0.866)
\qline(0.5,0.866)(1.25,-0.433)
\qline(1.25,-0.433)(0,0)
\qline(-0.25,-0.433)(1,0)
\qline(0,0)(1,0)
\psdot(0.5,-0.176) 
\qline(0.5,0.866)(0.5,-0.176)
\pscircle[fillstyle=solid,fillcolor=white](0.5,0){0.04}
\put(-0.2,-0.03){$A$}
\put(1.06,-0.03){$B$}
\put(-0.47,-0.48){$D$}
\put(1.3,-0.48){$E$}
\put(0.57,0.83){$C$}
\put(0.44,-0.38){$F$}
\put(0.42,0.08){$M$}
\psset{unit=3cm}
\endpspicture}

\def\AlternateInteriorAnglesFigure{%
\pspolygon[fillstyle=solid,fillcolor=yellow](0.7,0.7)(1.3,0.7)(0.8,0.425)
\pspolygon[fillstyle=solid,fillcolor=yellow](0.9,0.15)(0.3,0.15)(0.8,0.425)
\pspicture(3.2, 0.75)
\qline(0.0,0.15)(2,0.15)
\qline(0.0,0.7)(2,0.7)
\psdot(0.7,0.7)
\put(0.67,0.78){$p$}
\psdot(0.9,0.15)
\put(0.88,0.04) {$q$}
\psdot(0.3,0.15)
\put(0.3,0.04){$s$}
\psdot(0.8,0.425)
\put(0.85,0.36){$t$}
\qline(0.3,0.15)(1.3,0.7)
\qline(0.7,0.7)(0.9,0.15)  
\psdot(1.3,0.7)
\put(1.3,0.78) {$r$}
\put(-0.15,0.12) {$L$}
\put(-0.15,0.68) {$K$}
\endpspicture}

\def\FigureOneThirtyFiveColored{%
\psset{unit=1cm}
\pspicture(0,-0.5)(4,2.3)
\pspolygon[fillstyle=solid,fillcolor=yellow,linestyle=none]
(0,2)(0.5,0)(2.5,2)
\pspolygon[fillstyle=solid,fillcolor=yellow,linestyle=none]
(2,2)(4.5,2)(2.5,0)
\pspolygon[fillstyle=solid,fillcolor=red,linestyle=none]
(2,2)(2.5,2)(2.1,1.6)
\psdot(0.5,0)
\put(0.4,-0.35){$B$}
\psdot(0,2)
\put(-0.1,2.15){$A$}
\psdot(2,2)
\put(1.85,2.15){$D$}
\psdot(2.5,2)
\put(2.36,2.15){$E$}
\psdot(4.5,2)
\put(4.3,2.15){$F$}
\psdot(2.5,0)
\put(2.3,-0.35){$C$}
\qline(0,2)(4.5,2)  
\qline(0,2)(0.5,0)  
\qline(0.5,0)(2.5,2)  
\qline(2,2)(2.5,0)  
\qline(4.5,2)(2.5,0) 
\qline(0.5,0)(2.5,0) 
\psdot(2.1,1.6)
\put(2.24,1.47){$G$}
\endpspicture
}

\def\ProclusPostulateFourFigure{%
\pspicture(2, 0.8)(0,-0.06)
\qline(0,0.1)(1,0.1)
\qline(0.5,0.1)(0.5,0.75)
\qline(0,0)(1,0.2)
\qline(1.5,0.1)(2,0.1)
\qline(1.5,0.1)(1.5,0.75)
\put(-0.15,0.1){$H$}
\put(-0.15,-0.05){$K$}
\put(0.45,-0.02){$B$}
\put(0.95,-0.02){$C$}
\put(0.38,0.7){$A$}
\put(1.37,0.7){$D$}
\put(1.42,-0.02){$E$}
\put(1.92,-0.02){$F$}
\put(0.95, 0.23){$G$}
\endpspicture}

\def\OneSixteenInnerFigure{%
\pspicture(0,-0.05)(2,0.75)
\pspolygon[fillstyle=solid,fillcolor=yellow](1.8,0)(0,0)(2,0.8)(1.4,0.2)(1.8,0)
\pspolygon[fillstyle=solid,fillcolor=lightblue](0,0)(2,0.8)(1.2,0)(0,0)
\psline[linecolor=red](1,0.4)(1.8,0)  
\psline(0,0)(1.8,0)  
\qline(0,0)(0.8,0.8)  
\qline(0.8,0.8)(1.2,0)  
\psdot(0.8,0.8)   
\put(0.65,0.77){$A$}
\psdot(1.2,0)  
\put(1.17,-0.13){$C$}
\psdot(1.8,0)  
\put(1.75,-0.13){$D$}
\psline[linecolor=blue](1.2,0)(2 ,0.8) 
\psline[linecolor=blue](0,0)(1.19,0.477)  
\psline[linecolor=blue](1.21,0.483)(2,0.8)  
\psdot(1,0.4)  
\put(0.83,0.4){$E$}
\psdot(2,0.8)  
\put(2,0.67){$F$}
\pscircle[fillstyle=solid,fillcolor=white](1.4,0.2){0.03} 
\put(1.5,0.19){$H$}
\psdot(0,0)  
\put(-0.04,-0.13){$B$}
\endpspicture}

\def\OneSixteenOuterFigure{%
\pspicture(0,-0.15)(2,0.85)
\pspolygon[fillstyle=solid,fillcolor=yellow](0.8,0.8)(1.8,0)(0,0)(1,0.4)(0.8,0.8)
\pspolygon[fillstyle=solid,fillcolor=lightblue](0.8,0.8)(1.8,0)(1.2,0)(0.8,0.8)
\psline[linecolor=red](0.8,0.8)(1.17,0.50)  
\psline[linecolor=red](1.23,0.46)(1.8,0)  
\psline(0,0)(1.8,0)  
\qline(0,0)(0.8,0.8)  
\qline(0.8,0.8)(1.2,0)  
\psdot(0.8,0.8)   
\put(0.65,0.77){$A$}
\psdot(1.2,0)  
\put(1.17,-0.13){$C$}
\psdot(1.8,0)  
\put(1.75,-0.13){$D$}
\psline[linecolor=blue](1.2,0)(2 ,0.8) 
\psline[linecolor=blue](0,0)(1.19,0.477)  
\psline[linecolor=blue](1.21,0.483)(2,0.8)  
\psdot(1,0.4)  
\put(0.83,0.4){$E$}
\pscircle[fillstyle=solid,fillcolor=white](1.2,0.48){0.03}
\put(1.17,0.55){$J$} 
\psdot(2,0.8)
\put(2,0.67){$F$}
\psdot(0,0)  
\put(-0.04,-0.13){$B$}
\endpspicture}

\def\OneSixteenFigureTwo{%
\pspicture(0,-0.15)(2,0.75)
\pspolygon[fillstyle=solid,fillcolor=yellow](0.8,0.8)(1.2,0)(2,0.8)(1.2,0.48)(0.8,0.8)
\pspolygon[fillstyle=solid,fillcolor=lightblue](1.2,0)(2,0.8)(1,0.4)
\psline[linecolor=red](0.8,0.8)(1.8,0)  
\psline(0,0)(1.8,0)  
\qline(0,0)(0.8,0.8)  
\qline(0.8,0.8)(1.2,0)  
\psdot(0.8,0.8)   
\put(0.65,0.77){$A$}
\psdot(1.2,0)  
\put(1.17,-0.13){$C$}
\psdot(1.8,0)  
\put(1.75,-0.13){$D$}
\psline[linecolor=blue](1.2,0)(2 ,0.8) 
\psline[linecolor=blue](0,0)(1.19,0.477)  
\psline[linecolor=blue](1.21,0.483)(2,0.8)  
\psdot(1,0.4)  
\put(0.83,0.4){$E$}
\psdot(1.2,0.48) 
\put(1.17,0.55){$J$} 
\psdot(2,0.8)
\put(2,0.67){$F$}
\pscircle[fillstyle=solid,fillcolor=white](1.47,0.27){0.03}
\put(1.54,0.24){$H$}
\psdot(0,0)  
\put(-0.04,-0.13){$B$}
\endpspicture}

\def\FigureParallelDefFour{%
\psset{unit=2cm}
\pspicture(2,1.4)
\psline[linecolor=red](0.42,1.185)(0.935,0.715) 
\psline[linecolor=red](0.965,0.685)(1.485,0.215) 
\psline[linecolor=red](0.52,0.22)(0.935,0.685) 
\psline[linecolor=red](0.965,0.71)(1.4,1.2)  
\qline(0,0.2)(0.47,0.2)  
\qline(0.53,0.2)(1.47,0.2) 
\qline(1.53,0.2)(2,0.2) 
\qline(0,1.2)(0.37,1.2) 
\qline(0.43,1.2)(1.37,1.2) 
\qline(1.43,1.2)(2,1.2)  
\psdot(0,0.2)  
\put(-0.1,0.02){$C$}
\pscircle[fillstyle=solid,fillcolor=white](0.5,0.2){0.03} 
\put(0.45,0.02){$P$}
\pscircle[fillstyle=solid,fillcolor=white](1.5,0.2){0.03} 
\put(1.45,0.02){$Q$}
\pscircle[fillstyle=solid,fillcolor=white](0.4,1.2){0.03}  
\put(0.35,1.27){$R$}
\pscircle[fillstyle=solid,fillcolor=white](1.4,1.2){0.03}  
\put(1.35,1.27){$S$} 
\psdot(0,1.2)  
\put(-0.1,1.27){$A$}
\psdot(2,0.2) 
\put(1.95,0.02){$D$}
\psdot(2,1.2) 
\put(1.95,1.27){$B$}
\pscircle[fillstyle=solid,fillcolor=white](0.95,0.7){0.03} 
\put(0.88,0.8){$H$}
\endpspicture
\psset{unit=3cm}}

\def\ConnectivityFigure{%
\pspicture(2,0.6)
\psdot(0,0.25)
\psdot(2,0.25)
\psdot(1,0)
\psdot(1,0.5)
\qline(0,0.25)(0.2,0.25)
\psbezier(0.2,0.25)(0.4,0.25)(0.7,0.5)(1,0.5)
\psbezier(1,0.5)(1.3,0.5)(1.6,0.25)(1.8,0.25)
\qline(1.8,0.25)(2,0.25)
\psbezier(0.2,0.25)(0.4,0.25)(0.7,0)(1,0)
\psbezier(1,0)(1.3,0)(1.6,0.25)(1.8,0.25)
\put(-0.05,0.12){$A$}
\put(0.95,0.55){$B$}
\put(0.95,0.05){$C$}
\put(2.02,0.17){$D$}
\endpspicture}

\def\PizzaFigure{%
\pspicture(0.6,0.15)(0,-0.2)
\SpecialCoor
\pswedge[fillstyle=solid,fillcolor=lightblue,linestyle=none](0,0){0.3}{-30}{0}
\pswedge[fillstyle=solid,fillcolor=lightblue,linestyle=none](0,0){0.3}{15}{45}
\pswedge[fillstyle=solid,fillcolor=red,linestyle=none](0,0){0.3}{0}{15}
\pswedge[fillstyle=solid,fillcolor=red,linestyle=none](0,0){0.3}{-45}{-30}
\psarc(0,0){0.3}{-45}{45}
\qline(0,0)(0.3;0)
\qline(0,0)(0.3;15)
\qline(0,0)(0.3;45)
\qline(0,0)(0.3;-30)
\qline(0,0)(0.3;-45)
\endpspicture}

\def\PizzaFigureTwo{%
\pspicture(0.6,0.15)(0,-0.2)
\SpecialCoor
\pswedge[fillstyle=solid,fillcolor=lightblue,linestyle=none](0,0){0.3}{-15}{45}
\pswedge[fillstyle=solid,fillcolor=red,linestyle=none](0,0){0.3}{-30}{-15}
\pswedge[fillstyle=solid,fillcolor=red,linestyle=none](0,0){0.3}{-45}{-30}
\psarc(0,0){0.3}{-45}{45}
\qline(0,0)(0.3;45)
\qline(0,0)(0.3;15)
\qline(0,0)(0.3;-15)
\qline(0,0)(0.3;-30)
\qline(0,0)(0.3;-45)
\endpspicture}

\def\AngleEqualityFigureOne{%
\pspicture(0.5,0.4)(-0.5,0)
\SpecialCoor
\pspolygon[fillstyle=solid,fillcolor=lightblue,linestyle=none](0,0)(0.3;60)(0.3;120)
\qline(0,0)(0.5;60)
\qline(0,0)(0.5;120)
\psdot(0.5;60)
\psdot(0.5;120)
\psdot(0.3;60)
\psdot(0.3;120)
\psdot(0,0)
\qline(0.3;60)(0.3;120)
\put(-0.38,0.4){$A$}
\put(0.29,0.4){$C$}
\put(-0.13,-0.05){$B$}
\put(-0.27,0.2){$U$}
\put(0.18,0.2){$V$}
\endpspicture}

\def\AngleEqualityFigureTwo{%
\pspicture(0.5,0.4)(-0.5,0)
\SpecialCoor
\pspolygon[fillstyle=solid,fillcolor=lightblue,linestyle=none](0,0)(0.3;60)(0.3;120)
\qline(0,0)(0.4;60)
\qline(0,0)(0.4;120)
\psdot(0.4;60)
\psdot(0.4;120)
\psdot(0.3;60)
\psdot(0.3;120)
\psdot(0,0)
\qline(0.3;60)(0.3;120)
\qline(0.3;60)(0.3;120)
\put(-0.3,0.32){$a$}
\put(0.24,0.32){$c$}
\put(-0.09,-0.05){$b$}
\put(-0.27,0.2){$u$}
\put(0.18,0.2){$v$}
\endpspicture}

\def\AngleOrderFigure{%
\pspicture(0.5,0.4)(-0.5,0)
\SpecialCoor
\pspolygon[fillstyle=solid,fillcolor=lightblue,linestyle=none](0,0)(0.3;60)(0.3;120)
\qline(0,0)(0.45;60)
\qline(0,0)(0.5;120)
\qline(0,0)(0.75;30)
\psdot(0.5;30)
\psdot(0.46;60)
\psdot(0.75;30)
\psdot(0.5;120)
\psdot(0.4;120)
\psdot(0,0)
\qline(0.3;60)(0.3;120)
\put(-0.38,0.4){$J$}
\put(0.67,0.32){$K$}
\put(0.29,0.4){$X$}
\put(-0.33,0.3){$P$}
\put(-0.13,-0.05){$Q$}
\put(0.47,0.17){$R$}
\psline[linecolor=red](0.5;120)(0.75;30)
\endpspicture}

\def\OneTwentyFigure{%
\pspicture(-0.2,-0.2)(0.5,0.5)
\SpecialCoor
\pspolygon[fillstyle=solid,fillcolor=lightblue,linestyle=none](0,0)(0.5;-20)(0.25;240)
\qline(0,0)(0.5;60)
\qline(0,0)(0.5;-20)
\qline(0.5;60)(0.5;-20)
\qline(0,0)(0.25;240)
\qline(0.5;-20)(0.25;240)
\psdot(0,0)
\psdot(0.5;-20)
\psdot(0.25;240)
\psdot(0.5;60)
\put(-0.15,0){$A$}
\put(-0.28,-0.23){$B$}
\put(0.3,0.4){$D$}
\put(0.5,-0.25){$C$}
\endpspicture}

\title{Euclid after Computer Proof-Checking}   
\author{Michael Beeson}       
\date{\today
}  

\begin{document}

\maketitle

\begin{abstract}Euclid pioneered the concept of a mathematical
theory developed from axioms by a series of justified proof steps.
From the outset there were critics and improvers.  In this century
the use of computers to check proofs for correctness sets a new 
standard of rigor.  How does Euclid stand up under such an examination? 
And what does the exercise have to teach us about geometry, 
mathematical foundations, and the relation of logic to truth?  
\end{abstract}

\section{Introduction.} In approximately 300 BCE,  
a research institute known as the {\em Museum} 
 was founded at Alexandria.   A community of scholars lived 
and worked at the Museum, holding property in common and having
a common dining room.%
\footnote{According to Strabo's description 
(\cite{strabo}, Book~XVII,Chapter~ I).
See also see \cite[pp.~30,38]{sarton2}. }
Euclid joined this community,
 and became the author of several
influential books.%
The most famous of these is his
 {\em Elements}.%
\footnote{The story of how the {\em Elements} survived 
until modern times is fascinating but not strictly
relevant to our present topic; for the Library of Alexandria see  \cite{charlesriver,el-abbadi2020,
sarton2, strabo}; for its travels to India,
see \cite[Chapters 7 and 8]{oleary} and \cite[Book~XVII, Chapter 8]{strabo}.
For the Arabic flowering and return to Europe, see  \cite{oleary} and \cite[p.~595]{deRisi2016b}.
 For its printing in 1482 see \cite[p.~364]{heath-vol1}.
}

No other scientific book  has had 
an equal influence on the world. 
 The authors of the United States Declaration of Independence had
studied Euclid:  hence the ``self--evident truths'' from which 
the principles of the Declaration were deduced in Euclidean style.
Abraham Lincoln took several months studying Euclid%
\footnote{See \cite[pp.~64-65]{ketcham} and \cite[Chapter 14]{hirsch}.}
and thus learned what it means to demonstrate a proposition in court.
He put it to good use in the Lincoln--Douglas debates.  
Bertrand Russell studied Euclid at age eleven and ``did not know there
was anything so delicious in the world.''%
\footnote{See \cite{russell-autobiography} (Prologue).  But by 1902, when Russell was thirty, he
had changed his mind, and published an article pointing out errors
in Euclid~ I.1, I.4, I.7, and I.16 \cite{russell1902}.}  

From at least 450 CE (and probably before that!) there was no shortage
of critics, and suggestions for repair and improvement.  Dozens of 
these have been surveyed by De Risi \cite{deRisi2016b}.  Proclus already
criticized the parallel postulate (Euclid 5), saying that it did not 
deserve to be a postulate but should be proved, since we know that some 
kinds of (curved) lines can approach each other without meeting, so
why can't straight lines do that too?    In the nineteenth
century there was an increasing focus on the problem of whether the 
parallel postulate (Euclid 5) could be eliminated, by proving it from 
the other axioms.  Several famous mathematicians thought they had 
done so; Legendre published three different mistaken proofs.  This 
effort focused attention on rigor and careful deduction,  and 
led to careful axiomatic developments.  At the same time, logic 
was also being developed;  Boole's {\em Laws of Thought} was
published in 1853.   A milestone was reached with Pasch's work 
in 1882 \cite{pasch1882}, which introduced the notion of ``betweenness'' missing
in Euclid, and the idea that one should justify the existence of 
points where lines cross lines or circles, or circles cross other circles.

The geometers eventually realized that there is such a thing as 
non-Euclidean geometry, in which Euclid 5 fails but the other postulates 
hold.  This was the world's first independence proof.  The efforts
to codify systems of deduction,  stimulated by the pressing need to 
ensure correctness in geometry,  led directly to the logical systems
of Peano
 and Frege at the end of the nineteenth century, 
 which are the wellsprings of modern logic.  Geometry was the 
midwife of logic.

In 1899,  Hilbert published his very influential book \cite{hilbert1899}, in which he attempted to 
bring new standards of rigor to geometry.  
Hilbert's system mixed first-order and second-order
axioms with a little set theory thrown in for good measure,  but 
he clearly understood that axioms could be dependent or independent 
and a given axiom system might have different models.  Hence his 
famous dictum ``tables, chairs, and beer mugs,''  the point of which was 
 that reasoning should 
be fundamentally syntactic,  so that the steps could be checked independently
of the meaning of the terms; hence if you substitute ``tables, chairs, and 
beer mugs'' everywhere for ``points, lines, and planes,'' then everything 
should still be correct.  

Hilbert's book was about plane geometry,  but he did not follow 
Euclid.  His aim (as stated in the Introduction to \cite{hilbert1899})
was to ``establish for geometry a complete, and as simple as possible,
set of axioms.''  His method to achieve that goal was to show that 
arithmetic (addition and multiplication) can be defined geometrically.
This was first done more than two centuries earlier, in Descartes's famous
{\em La G\'eom\'etrie}, but Descartes's work was not based on 
specific axioms.   The Greeks never tried to multiply two 
lines to get another line;  the product of two lines was a rectangle.
Descartes broke through this conceptual barrier,  and Hilbert
made Descartes's arguments rigorous, using the theorems of the 
nineteenth-century jewel ``projective geometry.'' 
Hilbert's book put in the mathematical bank the profits of the 
discoveries about non-Euclidean geometry: geometry was no longer about 
discovering the truth about points, lines, and planes,  but instead it 
was about what theorems follow from what axioms.  Rather than deduce
Euclid's theorems directly from his axioms, Hilbert relied on the 
indirect argument that every geometrical theorem could be derived by 
analytic geometry, and since coordinates could be introduced by 
the geometric definition of arithmetic, every theorem had a proof from 
his axioms.   

In 1926--27, Tarski lectured on his first-order theory of geometry.%
\footnote{According to \cite[p.~175]{tarski-givant}. This reached
journal publication only in \cite{tarski1959}. }
But Tarski and his collaborators too ignored Euclid, focusing as Hilbert did on the development of geometric arithmetic and the characterization of models.

Less than a decade after the first electronic computers were 
available, there were efforts to apply them to the problem of finding
proofs of geometry theorems (\cite{gelernter1959, gelernter1960}).  Later, the verification
of geometry theorems by algebraic calculation became an advanced art form.
(See the introduction to \cite{beeson2019} for citations and history).  On the other hand,  since the 1990s there has been a body of 
work in ``proof-checking,''  which means that a computer program 
checks that the reasoning of a given proof is correct (as opposed to 
``automated deduction,'' in which the computer is supposed to find the 
proof by itself).  In recent years there have been several high-profile
cases of proof-checking important theorems whose large proofs involved
many cases.  Tarski's work on geometry, as presented in \cite{schwabhauser},
has also been the subject of proof-checking and proof-finding experiments
\cite{beeson2017a,boutryphd,narboux2017b,narboux2017c,narboux2017}.  

Until 2017, nobody had tried to proof-check Euclid directly.  
That omission was remedied in \cite{beeson2019}.
We produced formal proofs of 245 
propositions, including the 48 propositions of Euclid Book I, and 
numerous other propositions that were needed along the way, and 
 translated those proofs into the languages of 
two famous proof-checkers, Coq and HOL Light.  
The purpose of this article is to consider what is to be learned from 
that proof-checking experiment.  That work invites the
question, {\em How wrong was Euclid?}   That 
question will be answered below.  Here we consider as well:

\begin{itemize}
\item  {\em What was Euclid thinking?} 
\item  {\em What is the relation between logic and truth?}
\item  {\em How should we choose the primitive notions for a formalization?}
\item  {\em What level of mathematical infrastructure should be visible?}
\end{itemize}

Relevant to those questions are the following principles of 
mathematical practice (which are sometimes followed and sometimes not):

\begin{itemize} 
\item {\em If it can be defined, instead of taken primitive, it should be defined.}
\item {\em If it can be proved, instead of assumed, it should be proved.}
\end{itemize}

Tarski clearly believed in both these principles: he famously sought to 
minimize the number of primitives and the number of axioms, and 
tried to find axioms that could be stated in the primitive language without
preliminary definitions.  Hilbert evidently wasn't such a strong 
believer; his axioms were stronger by far than they had to be,
because
he wanted to reach his goal of defining arithmetic as quickly as possible.
Of course, it was the second of those principles that drove people to 
work so hard trying to prove Euclid~5 from the other axioms, and 
over the centuries there were many claims that Euclid~4  (all right 
angles are equal) can be proved (this will be discussed below).

I have come to think of the proof-checking of Euclid as the 
enterprise of {\em providing infrastructure}.   Euclid's proofs
are like mountain railroads:  there is a clear origin and destination,
but there are a lot of possibilities to go off the rails.  There are
places in Euclid where the train would fall into the abyss;  there are
other places where the trestles are very weak.  After \cite{beeson2019},
we have 245 propositions to support Euclid's 48.   It is a much 
sturdier system.   Here I will discuss it and consider its 
relation to Euclid's {\em Elements}, Book~I, which it is meant to formalize.

\section{Structure of Euclid's {\em Elements}.}
Euclid had {\em definitions}, {\em common notions}, {\em axioms}, and 
{\em postulates}.
Nowadays, the  common notions,
axioms, and postulates would be lumped together and considered axioms.
In Euclid, the common notions were intended to be principles of 
reasoning that applied more generally than just to geometry.  
For example, what 
we would now call equality axioms,  or such principles as 
``the part is not equal to the whole'' or ``the whole is equal to 
the sum of the parts.''  Any kind of thing might have parts, which were
always the same kind of thing: the parts of a line were smaller lines, etc.
 The axioms and postulates were about geometry.
The distinction between an ``axiom'' and a ``postulate,'' according 
to Proclus \cite[\S 201, p.~157]{proclus},  who attributes the distinction
to the earlier mathematician Geminus, is that a postulate asserts that some point can be constructed,
while an axiom does not.  In modern terms, an ``axiom'' is purely universal,
while a postulate has an existential quantifier.  Euclid's Postulate~4
(all right angles are equal) had no associated construction, and for that
reason, critics in antiquity (such as Proclus, in \cite[\S 192, p.~150]{proclus} said it did not deserve
the name ``postulate.''

Heath's translation lists five common notions, five postulates, and zero
axioms.  Simson's translation \cite{simson} lists three postulates, twelve axioms,
and zero common notions.  The extra axioms are discussed by Heath \cite[p.~223]{euclid1956},
where they are rejected.  During the centuries since Euclid, a great many
axioms have been put forward as ``should have been included.''  The 
axiomatization of \cite{beeson2019}, 
discussed in this article,  is thus the last elephant in a long parade.%
\footnote{See \cite{deRisi2016b} for a thorough discussion.  In that
paper, it takes De Risi twenty pages just to list the different axioms
that various authors have added to Euclid's. 
}
  But the 
seal of correctness given by computer proof-checking does support the  
claim to be the ``last elephant.''   That axiomatization followed the consensus of the 
centuries on these points:  Postulate~4 should be proved,  the SAS criterion of 
Proposition~I.4 should be an axiom,  and Postulate~5 has to remain.
It also followed the consensus of the twentieth century that the concept 
of betweenness is needed, and  it adopted axioms similar to those of 
Hilbert and Tarski about betweenness.

\section{Lines.} 
Euclid's lines were all finite lines; a line could be extended, but 
after the extension it was still a finite line, just longer.%
\footnote{Nowadays lines are infinite; Euclid's (finite) lines are
today's ``segments.''  A ``line'' for Euclid could also be curved; hence
``straight line'' adds something.  However, the default in Euclid, though 
not necessarily in all Greek geometry, is that lines are straight. 
We follow Euclid.} 
  How far
could it be extended?  Later Hilbert and Tarski both answered that question 
by taking it as an axiom that $AB$ could be extended by an amount equal to a given 
line $CD$.   This is sometimes known as the {\em rigid compass}:  we can 
measure off $CD$ with a compass, then move the compass to the end of $AB$
and lay off the measured amount along a straightedge.   By contrast, 
Euclid used a {\em collapsible compass}:  you could only extend a line
by putting one tip of the compass on $B$ and the other tip on some 
existing point, then drawing a circle; then $AB$ would be 
``drawn through'' to meet that circle.  Euclid's second proposition, 
I.2,  gives a beautiful proof that the rigid compass can be simulated
by a collapsible compass;  that conclusion is axiomatized away by 
Hilbert and Tarski, violating the principle ``if it can be proved, it 
should be proved.''%
\footnote{A referee pointed out \cite{rigby1970}, which gives an axiomatic
framework for a part of geometry in which the collapsible compass is preserved.}  

In Euclid's time it had not yet been realized that
not everything can be defined:  certain notions must be taken as primitive,
because otherwise there will be nothing to use in the first definition.
Euclid defined a line as that which has length but no breadth.  Other 
sources make clear that the Greeks regarded lines (curved or straight) as
the traces of a motion.%
\footnote{This view of lines goes back at least to 
Aristotle; see \cite[p.~79]{proclus}, where Proclus says Aristotle
regarded a line as ``the flowing of a point.''}
   In particular, most modern mathematicians, 
steeped in set theory from their youth,  think that a line is equal to 
the set of its points.  This was definitely not the Greek conception.
This viewpoint was regarded as thoroughly refuted by the paradoxes 
of Zeno, which were already centuries old by the time of Euclid.  
These paradoxes depended on the conception that, if a line were made up
of its points, the points would be like tiny beads on a string.  

In mathematics there is another principle: {\em Perhaps you do not need to 
know what it is,  if you know how to use it.}  Euclid never once appeals
to his definition of a line, so in essence he did treat it as a primitive 
notion.  He always refers to a line by two distinct points, as in ``line $AB$.'' 
You never see ``line $\ell$\,'' in Euclid.   Lines are used to construct other 
points, as the intersection points of lines with other lines or with circles.

Euclid spoke of two lines being ``equal.''   He meant by this what 
Hilbert called ``congruent.''   Probably he had in mind that a rigid motion 
could move one line to coincide with the other.   It seems certain that 
he did not mean that they had the same length, as measured by a number.   

The other obvious primitive notion about lines (besides equality)
 is  the notion ``point $P$ lies on line $AB$.'' 
This is the notion of {\em betweenness}, which was explicitly axiomatized
by Pasch only in 1882.  (One can either use {\em strict} betweenness, as 
Hilbert did, or nonstrict betweenness, as Tarski did.)   Both Hilbert and 
Tarski had some axioms about betweenness.  Euclid, on the other hand, 
had none.   

Euclid did use the notion that one line is ``less than'' or ``greater than''
another.   That notion can be defined in terms of betweenness and equality:
$AB < CD$ if there is some $X$ between $C$ and $D$ such that $AB$ is equal 
to $CX$.  Conversely, we could define (for collinear points)
 $X$ is between $C$ and $D$ if $CX < CD$ and $XD < CD$.
For Euclid, ``less than'' was a natural notion, because $AB < CD$ meant
that $AB$ was equal to a part of $CD$.  The notion of one thing being 
part of another thing was regarded as a ``common notion,'' i.e., a notion 
that applied to things in general, as opposed to notions specific to geometry.
Euclid took ``part'' as a primitive notion, making no attempt to define it, 
but it was clear that the parts of a thing were the same kind of thing: the parts of a line were lines,  the parts of an angle were angles.  The definition of 
{\em point} is {\em that which has no parts}.   The common notions probably 
seemed completely precise to Euclid, because he felt that the notion of ``part''
was clear.  From the modern point of view, there are questions:  if $AB$ is 
divided into two parts by its midpoint $M$,  what about point $M$ itself?  
Is it in both parts or neither part?  Is the whole $AB$ really the ``sum'' of 
$AM$ and $MB$?  This problem comes back ``in spades'' when we divide a 
triangle into two parts:  what about the separating line?  When we put 
the parts back together, there would be no cut remaining.

We chose to follow Pasch, Hilbert, and Tarski in using betweenness as a 
primitive notion; we followed Hilbert in using strict betweenness, because
Euclid has no ``null lines'':  when he says ``line $AB$,'' $A$ and $B$ are
always distinct points.  Of course, nonstrict betweenness is easily 
defined in terms of strict betweenness, so it is really inconsequential which 
is taken as primitive.  We write $\B abc$ for ``$b$ is between $a$ and $c$.''
The three betweenness axioms and their names are given here:

\begin{center}
\begin{tabular}{l l}
symmetry & $\B abc \leftrightarrow \B cba$ \\
identity & $\neg\, \B aba$ \\
inner transitivity & $\B abd \ \land \ \B bcd \rightarrow \B abc$
\end{tabular}
\end{center}

It turns out that the theory of betweenness is not quite trivial.  There are
delicate and interesting questions about the axiomatization of this notion.
Tarski's original version of his theory had 16 axioms.  Working with 
his students,  in 1956--57 several of these axioms were proved
dependent on others.  That the three axioms given above are enough is
remarkable, but it is not very important for formalizing Euclid, 
 since
Euclid never even mentioned betweenness.  We must add the required
infrastructure, but it won't matter exactly how we do it, so we do not
go further into the matter here.%
\footnote{See \cite{tarski-givant} for a complete history of Tarski's
axiom systems, including the discoveries of dependencies among the axioms.} 
If five betweenness axioms had been required instead of three, we would have 
just added five betweenness axioms.  

\begin{figure}[ht]
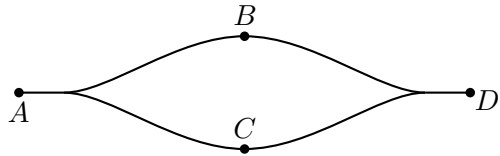

\caption{The connectivity axiom: two lines cannot enclose a space.}
\label{figure:connectivity}
\center{\ConnectivityFigure}
\end{figure}

The axiomatization of \cite{beeson2019} supplemented Euclid's
axioms, as given in the Heath translation, with the axiom the authors called {\em connectivity}:  if $B$ and $C$ are 
both between $A$ and $D$,  but neither $\B BCD$ nor $\B CBD$,
then $B=C$.  In other words, the picture in Figure~\ref{figure:connectivity}
is impossible.  This axiom is found in several of the Greek manuscripts
used by Heiberg in preparing his influential translation.%
\footnote{See \cite[p.~638]{deRisi2016b}. One of these manuscripts
is the oldest copy of Euclid to survive into the modern age.} 
  The ancients
did not use betweenness; they expressed this axiom as ``two lines intersect in at most one point,'' or as ``two lines cannot enclose space.''  The latter
lacks precision as ``space'' is not defined; the former is equivalent to 
our axiom.  Zeno of Sidon attacked Proposition~I.1 on the grounds that 
it is not conclusive unless it first be assumed that neither two straight lines
nor two circumferences can have a common segment.%
\footnote{See \cite[p.~359]{heath-vol1}, who in turn cites Proclus.} 

The connectivity axiom is used to prove the uniqueness of the midpoint, 
and the basic property of collinearity, if $a$, $b$, and $c$ are collinear
and $a$,$b$, and $d$ are collinear, 
and $b \neq c$ and $a \neq b$, then $b$, $c$, and $d$ are collinear.
That fact is used in almost
every one of our formal proofs. 
How did Euclid get by without that axiom?
By not proving the uniqueness of the midpoint, and using the properties
of collinearity without explicit mention.

If we had adhered to the principle, ``If it can be proved, it should 
be proved and not assumed,''  we would have formalized Gupta's 
1965 proof of the connectivity axiom,  which was proved in 
his thesis \cite{gupta1965} and can more easily be found as Satz~5.3 in
\cite{schwabhauser}.  But that proof, with its elaborate counterfactual
diagram, is not ``in the spirit of Euclid,'' so we chose to include
the axiom of connectivity, which so many over the centuries have thought
should be included.%
\footnote{Hilbert smuggled that axiom into his system by requiring
uniqueness in his angle-copying axiom.  Pasch had it explicitly
as an axiom, Kernsatz~V, \cite[p.~5]{pasch1882}. 
}  

Collinearity, which is defined from betweenness by enumerating the cases, 
plays a larger-than-life role in our computer formalization; mostly
in the form of noncollinearity.  There are many cases in which to 
state a theorem formally, one must add to the hypothesis statements
like ``$ABC$ is a triangle.''  Since we define a triangle as three
noncollinear points, this amounts to ``$A$, $B$, and $C$ are not collinear.''   
More than 
half the individual inferences in our formalization turned out to 
be statements of collinearity or noncollinearity.  These statements
are ``pure infrastructure,'' absolutely necessary to prevent the 
theorems from collapsing into the abyss,  but absent in Euclid, 
and serving only to ensure that the diagram does not degenerate.
Euclid simply assumed that points that appear to be on a line are 
indeed on that line, and points that appear noncollinear, are noncollinear.  

\section{Angles.} \label{section:angles}
The Greek concept of angle was more general than the modern concept,
as we see from Euclid's definitions of {\em planar angle} (the inclination
to one another of two lines in a plane that meet one another and do not 
lie in a straight line), from which we see that one could also consider
angles that do not lie in a plane; and of {\em rectilinear angle}, 
when the lines are straight lines.   In particular, two touching 
circles form an angle that is not rectilinear; modern mathematics does 
not use the word ``angle'' in that situation.  Neither, as it turns out,
does most of Euclid; and here we use ``angle'' for ``rectilineal angle.''
Euclid always refers to angles by three points, as in ``angle $ABC$,'' 
never using the more modern notation ``angle $\alpha$.''  As with lines,
he uses the word ``equal'' instead of Hilbert's ``congruent.''  

As with lines, Euclid never once appeals to the definition of ``angle,''
so we must ask how angles are used,  rather than what they are.  
It turns out that (in addition to the equality relation) there are
relations of ``less than''  and ``greater than'' between angles, 
and two angles can be ``taken together,'' or sometimes ``added,'' 
in such a way that the common notions ``if equals are added to equals
the results are equal,'' and ``the whole is equal to the sum of its parts''
are applicable to angles.  In effect, then,  Euclid treats angles
as a primitive notion, with an ordering relation.   One consequence of 
Euclid's definition is that all his angles are ``less than two right angles,''
or in modern terminology, less than 180 degrees.  For otherwise 
angle $ABC$ would not be determined, unless we were to insist that 
angle $ABC$ is not the same as angle $CBA$; but Euclid clearly considers
$ABC$ equal to $CBA$.  

Although Euclid refers to angles by the names of three points, sometimes
we have to consider different ways of choosing those three points.
To give a specific example, consider angles $BAE$ and $BAC$, where 
$E$ lies between $A$ and $C$. 
(See Figure~\ref{figure:OneSixteenOne}). In the proof of I.16, Euclid
implicitly equates those angles.  Did he consider them to be 
 the {\em same}
angle (identical), or  merely {\em equal} angles?  This question 
cannot be answered by reading Euclid, since he never explained what he 
meant by ``equality.''  In practice, it isn't going
to matter whether we consider them to be equal angles, 
or are two names for the same angle.  The same missing steps will have
to be supplied.   

Euclid never mentions ``rays,''  because all his lines are finite.  
Therefore Euclid's angles all have finite (but extensible)
sides.  Nevertheless the  language of ``rays'' is useful; we say 
$x$ lies on the ray $AB$ if $\B ABx$ or $x = B$ or $\B AxB$.   

\begin{figure}[ht]
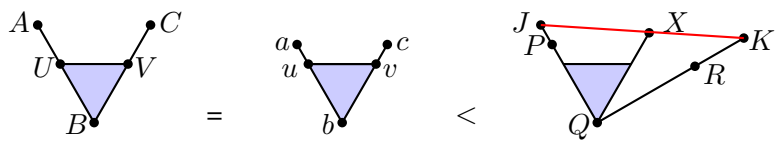

\caption{Definitions of angle equality and angle less than.}
\label{figure:angles}
\center{\AngleEqualityFigureOne = \AngleEqualityFigureTwo
  $<$ \AngleOrderFigure}
\end{figure}

At some point it was realized that angles could just be {\em defined}
as triples of points.%
\footnote{
As far as I can determine, Mollerup \cite{mollerup1904} deserves the credit 
for these definitions. 
}
Angles then are a case in point for the principle, ``if it can be defined,
it should be defined.''  In our computer proof-checking of Euclid, we 
used Tarski's  points-only language, so lines became pairs of points, 
and angles became triples of points.  
 Then the notions of 
equal angles and ``less than'' for angles can also be defined. 
See Figure~\ref{figure:angles}, in which the blue triangles are congruent. 

Angles $ABC$ and $abc$ are equal if
there exist points $U$, $V$, $u$, and $v$ on rays $BA$, $BC$, $ba$, and $bc$
respectively, such that $BU = bu$ and $BV = bv$ and $UV = uv$. 
And angle $abc$ is less than $PQR$ if there are $X,J,K$ with
$\B JXK$ and $J$ on ray $QP$ and $K$ on ray $QR$
such that $abc = PQX$.

If we view formalization as providing infrastructure, there is a lot
of infrastructure connected with equality and inequality of angles.
To start with, we must verify that angle equality is an equivalence
relation; that requires Euclid I.4, the SAS congruence criterion,
which is discussed in the next section.  This particular piece of 
infrastructure results from having to prove what Euclid took as 
``common notions''.  But not all the required angle infrastructure
is of that nature.  
For example, in the proof of Euclid~III.20 there is 
an unjustified step; I mean by ``unjustified''  that
Euclid did not write any justification for it, in the sense of a 
reference to an earlier proposition or definition.  To justify it,
he would have had to prove a proposition something like ``the sum 
of the doubles is the double of the sum,''  or more explicitly,
if two angles are doubled, then their doubles taken together equal
the double of the angles taken together.  In more modern terms, 
if we call the angles $1$ and $2$,  twice the sum of angles $1$ and $2$
is the sum of twice angle $1$ and twice angle 2.  
\begin{figure}[ht]
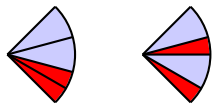

\caption{The sum of the doubles is the double of the sum.}
\label{figure:pizza}
\center{\PizzaFigureTwo  \PizzaFigure}
\end{figure}
If someone had pressed Euclid on this point, he would have justified this
step by ``the whole is equal to the sum of the parts.''   If we have 
four slices of pizza labeled 1,2,1,2  in order, and we take them 
out of the box and then put them back in the order 1,1,2,2,  lo and behold,
they fit into the original angle exactly. See Figure~\ref{figure:pizza}.  Now try to prove it 
by the methods of high school geometry.  That is ``infrastructure.''
You see from this example that it is logic and the choice of axioms that
give rise to the need for infrastructure.   Computer proof-checking
only shines a light on the situation.  You cannot convince the computer
by a story about pizza.

To formalize Euclid~III.20,  one has to define addition and subtraction
of angles;  Euclid thought that ``taken together'' was clear enough without
a definition.  It was a common notion, applying to any kind of ``thing,''
geometric or not.  This was not an ``operation'' in the modern sense, so 
the questions of commutativity and associativity were not considered.
But they now require proof.  Since Euclid had no notation for addition 
(of any kind of object), his notation was completely relational; that is,
he had to say $ABC$ and $abc$ taken together are equal to $PQR$, since
he could not say $ABC + abc$.  Algebraic (functional) notation for
addition and multiplication was introduced some 1800 years
after Euclid.  
To discuss the sum of four angles with different order and associativity
is extremely awkward in Euclid's notation, which was, of course, reflected
in our formal notation since we purposely tried to match our notation to 
Euclid's.  

Whether we take angles as primitive (as Hilbert did) or defined (as we do,
and Tarski did) doesn't matter very much for the formalization of Euclid,
as if one takes angles as defined, then one proves their basic properties and 
from then on things look pretty much the same as if angles were primitive. 
The exact choice of axiom system is not philosophically or mathematically important
(though it might make certain proofs easier or harder);  what is important is the complete
precision of detailed proofs.

\section{The SAS congruence theorem and the 5-segment axiom.}
Euclid attempted, in Proposition I.4,  to prove the side-angle-side
criterion for angle congruence (SAS).  But his ``proof'' appeals to 
the invariance of triangles under rigid motions, about which there is 
nothing in his axioms, so for centuries it has been recognized that 
in effect SAS is an axiom, not a theorem.

Before discussing SAS, we discuss the notion of triangle congruence.
Euclid does not define either ``triangle'' or ``congruent triangles,''
but taking I.4 to define the SAS criterion, the conclusion includes
the pairwise equality of corresponding sides,  and the pairwise equality
of corresponding angles.  As we discussed in Section~\ref{section:angles},
equality of angles is taken as a defined notion, rather than primitive.
With the definitions given there,  if two triangles have all pairs of 
corresponding sides equal, they automatically have corresponding angles equal too.
Hence SAS needs only to mention the equality of corresponding sides in 
its conclusion.  Moreover, two of those pairs of sides are equal by hypothesis,
so the conclusion of SAS just needs to be one equality of lines.  Next
we discuss how to formulate SAS without explicitly mentioning angles.

To do that, we use an axiom known as the ``five-line
axiom.''  This axiom is illustrated in Figure~\ref{figure:5-segment}.
The point of this axiom is to express SAS without mentioning angles at all.
To understand the relationship of SAS to the five-line axiom, let us
express Euclid I.4 (which is SAS) using Figure~\ref{figure:5-segment}. 
The hypothesis  is that $db = DB$, $dc=DC$, and  angles $dbc$ and $DBC$
in Figure~\ref{figure:5-segment} are equal.  The conclusion 
is that $dc = DC$.  The point of the 5-line axiom is to replace the 
hypothesis ``angle $dbc$ and $DBC$ are equal'' by the hypothesis that 
triangles $abd$ and 
$ABD$ are congruent, i.e., $ab=AB$, $bd=BD$, and $ad=AD$.   To rephrase the matter:
 The  hypothesis of the five-line axiom expresses the congruence (equality, in Euclid's
phrase) of angles $dbc$
and $DBC$ by means of the congruence of the exterior triangles $abd$ and 
$ABD$.  

It is not difficult to derive Euclid's I.4 from the five-line axiom.  It is 
also not difficult to derive the five-line axiom from I.4.  So,  it just 
a choice whether to take the five-line axiom, or I.4,  as an axiom, and 
after deriving one from the other, it makes little difference to the 
subsequent development.  We chose to follow Tarski in using the five-line 
axiom, since it can be stated succinctly using the primitive notions of the 
language (without abbreviations).  
  
\begin{figure}[ht]
\center{\TarskiFiveSegmentFigure}
\caption{If the four solid lines on the left are equal to the 
corresponding solid lines on the right, then the dashed lines
are also equal.}
\label{figure:5-segment}
\end{figure}
\FloatBarrier
This version of the five-line axiom was introduced by Tarski,
although we have changed nonstrict betweenness to strict betweenness.%
\footnote{The history of this axiom is as follows.
The key idea (replacing reasoning about angles by reasoning
about congruence of segments) was
introduced (in 1904) by J. Mollerup \cite{mollerup1904}.
His system has an axiom closely related to the 5-line axiom,
and easily proved equivalent.  Tarski's version \cite{tarski-givant}, however, is 
slightly simpler in formulation.   Mollerup (without comment)
gives a reference to 
Veronese \cite{veronese1891}.  Veronese does have a theorem 
(on page 241) with the same diagram as the 5-line axiom, and
closely related,  but he does not suggest an axiom related to this 
diagram.}

\section{Line-circle and circle-circle continuity.}
Euclid's  first proposition, Proposition~I.1,  
constructs an equilateral triangle by drawing two circles 
of radius $AB$ with centers at $A$ and $B$, respectively.  A meeting point
of these circles is a point $C$ equidistant from $A$ and $B$.  But 
why do the circles meet?  Euclid smuggles the point $C$ into the proof
by using a ``definite description'':  ``the point $C$ at which the circles
cut one another.''  The modern consensus is that this is a case of 
a ``missing axiom.''  We have to supply the {\em circle-circle continuity axiom},
according to which, if one circle has points inside and outside the other
circle, then the two circles meet.%
\footnote{De Risi \cite[p.~614]{deRisi2016b} gives credit to Richard
for first recognizing (in 1645) that an axiom would be required, though others
previously noted the gap in the proof of I.1, without filling it.} 
 The words {\em inside} and {\em outside}
are defined as follows: $X$ is inside a circle centered at $O$
if $OX$ is less than some 
radius, and outside if $OX$ is greater than some radius.  A radius is a line
connecting the center with a point on the circle. 

There is also the { line-circle continuity axiom},  asserting that if 
$A$ is inside a circle, and $P$ is any point different from $A$,
then there are two points collinear with $AP$ lying on the circle,
and one of them has $A$ between it and $P$.  This axiom is needed
twice in Euclid Book~I, in I.2 and I.12,  where the ``dropped perpendicular''
to a line from a point not on the line is constructed by drawing 
a sufficiently large circle,  which must meet the line in two points, 
forming a line whose perpendicular bisector is the desired perpendicular. 

These axioms are used seldom, but crucially, in Euclid.  Specifically,
circle-circle is used in Euclid~I.1,  which is the ``bootstrap'' proposition
for the first ten; it is used again in I.22, to construct a triangle out 
of three given lines.  And dropped perpendiculars, constructed by 
line-circle, are of course fundamental.

In modern times it has been shown that, using only the other axioms,
line-circle implies circle-circle and vice versa.   The only purely 
geometric proof known of these facts makes use of the ``radical axis''
and is a little complicated; see \cite{hartshorne}.  According to 
the principle ``if it can be proved, it should be proved,''  we should 
have just taken one of these axioms; but instead we took both
line-circle and circle-circle.  Had we taken only one, we would have
not been proof-checking Euclid, but proof-checking the modern theorems
about the radical axis.  That could clearly be done, but would not 
have added anything significant.

There is another axiomatic question about line-circle and circle-circle.
They were mentioned by Tarski in his original paper, but not included
in his full axiom list; he probably thought they followed from his 
Dedekind-style continuity schema (A11), which asserts that first-order
Dedekind cuts are filled.  This seems plausible until you actually try 
to prove it.  But to do so you need to drop perpendiculars to a line,
and to do that you need, if you follow Euclid,  circle-circle 
intersection.  So the argument is circular.   This is fixed by 
appealing to the 1965 thesis of Gupta \cite{gupta1965}, whose proofs 
were finally published in \cite{schwabhauser}.  Gupta showed,  amazingly, how
to construct both dropped and erected perpendiculars
 {\em without using circles at all}. So it does turn out to be correct
to omit circle-circle in favor of (A11), but Tarski certainly didn't
have a proof of that in 1927 or even 1959.

\section{Book Zero.}
It may surprise the reader when I say that even after the point $C$ in 
Proposition~I.1 has been admitted to exist, there is yet another defect
in Euclid's proof.  Namely, although we now know $AC=AB=BC$, in order
to prove $ABC$ is a triangle, we must also prove that $A$, $B$, and $C$
are not collinear.  To prove this we used a lemma we called
{\em partnotequalwhole}:  if $\B ABC$ then $AB \neq AC$.%
\footnote{The name of the lemma is taken from Euclid's Common Notion 5;
but Euclid does not cite CN5,  in I.1 or anywhere else, and De Risi
\cite{deRisi2020} after careful study reaches the conclusion that the original Euclid 
had only the first three of the five common notions given in the Heath
translation.  It seems, however, that he {\em should} 
have had CN5, and should have cited it in I.1.}
The proof
needs nine inferences, starting with Euclid's extension axiom to 
extend $ABC$ to another point $D$ left of $A$, i.e., with $\B DAB$.
  Then we need to 
show that $DABC$ occur in that order, in particular $\B DAC$,
using one of the lemmas about 
betweenness alluded to in the section on betweenness.   These
lemmas are part of what we might call ``Book Zero'';  Book Zero
contains infrastructure that is more fundamental than the propositions
of Book I.  It consists of the betweenness lemmas, several lemmas with 
a similar flavor to {\em partnotequalwhole},  and several trivial lemmas
that reflect the fact that we represented lines as pairs of points.
Thus if $AB=CD$, we also have $BA=CD$ and $AB=DC$ and $BA=DC$,
expressing the fact that these are unordered lines, not vectors.  

The fundamental properties of collinearity and noncollinearity, which 
are never mentioned in Euclid, should also be considered part of Book Zero.

\section{Pasch's axiom.}
Pasch \cite{pasch1882} not only introduced betweenness, but also 
the axiom that later was given his name.%
\footnote{Specifically, Kernsatz~IV, \cite[p.~20]{pasch1882}.
Pasch called his axioms ``Kerns\"atze.''  The ``kernel'' of a theory
consisted of kernel concepts and kernel theorems, but Pasch had a 
modern understanding of completeness and consistency, as p.~18 indicates.
} 

\begin{figure}[ht]
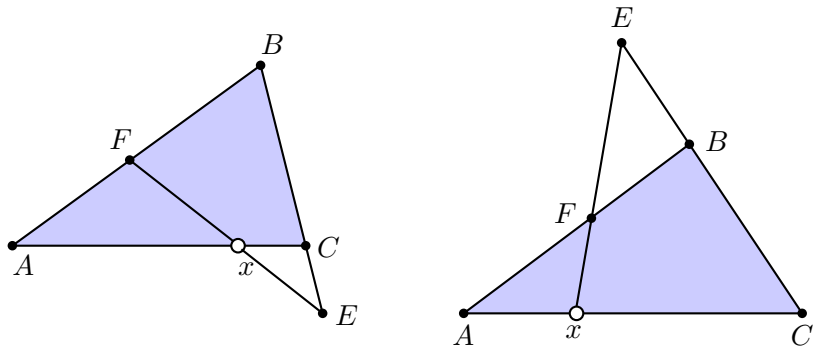
  
\caption{ Inner Pasch (left) and outer Pasch (right).  Line $EF$ meets triangle $ABC$ in one side $AB$, and meets an extension of side $BC$.  Then 
it also meets the third side $AC$.
 The open circles show the points asserted to exist. }
\label{figure:InnerOuterPaschFigure}
\InnerOuterPaschFigure
\end{figure}

  Pasch's axiom (Figure~\ref{figure:InnerOuterPaschFigure})
says that if $ABC$ is a triangle, and line $DE$ lies in a
plane with $ABC$ and meets $AB$
in a point $F$ between $D$ and $E$, then $DE$ or an extension of $DE$
meets $AC$ or $BC$.%
\footnote{In fact the axiom was formulated two centuries earlier
by Roberval; see \cite[p.~632, axiom C18, and discussion
p.~615]{deRisi2016b}.  One may wonder why it took two millenia
for this axiom to be formulated.  See \cite{deRisi2019} for a 
penetrating historical and philosophical discussion of that question.
In Pasch's statement, the first ``between'' here used ``innerhalb''
and so was strict; the second did not use ``innerhalb'' so was not
strict betweenness.  Pasch used both.)
}

Pasch's requirement that $DE$ lie in a plane with $ABC$  of
course cannot be dropped, since the line might not lie in the plane 
of the triangle.  In order to drop that hypothesis, obtaining
a statement that mentions only betweenness, one must strengthen the hypothesis
so that the line certainly lies in the plane of the triangle. There
are two ways to do this, resulting in axioms known as ``outer Pasch''
and ``inner Pasch.''  See Figure~\ref{figure:InnerOuterPaschFigure}.
In Pasch's own axiom (and figure) there is no requirement for point $E$
to be collinear with $BC$; that was added by Peano to make the coplanarity
hypothesis explicit.

These planar forms of Pasch's axiom were invented by Peano and published
in 1889,  seven years after Pasch.%
\footnote{
Axiom XIII in \cite{peano1889} is outer Pasch, with $\B abc$
written as $b \in ac$.  
Axiom XIV is inner Pasch.  Peano wrote everything in formal symbols 
only, and eventually bought his own printing press to print his books himself. See \cite{kennedy2002}.}

Inner Pasch has a certain symmetry: In Figure~\ref{figure:InnerOuterPaschFigure},
we could just as well have shaded triangle $BEF$ instead of triangle
$BAC$.  One soon becomes accustomed to noticing (and shading) the whole
quadrilateral.  

Gupta's thesis (which contained enough material for three theses) 
contains proofs that outer Pasch implies inner Pasch, and vice versa,
using the other axioms (but not continuity).   Of course Euclid,
who never mentions betweenness, did not explicitly use either version of 
Pasch, but they both came up naturally when we formalized Euclid.  
As with the continuity axioms, we could have picked just one, but then 
we would have been proof-checking Gupta's thesis, in addition to Euclid;
so we just took both inner and outer Pasch as axioms.
 The need for the Pasch axioms is pervasive:  we used 
 inner Pasch 36 times and outer Pasch 31 times in formalizing
 Euclid Book~I.

\begin{figure}[ht]
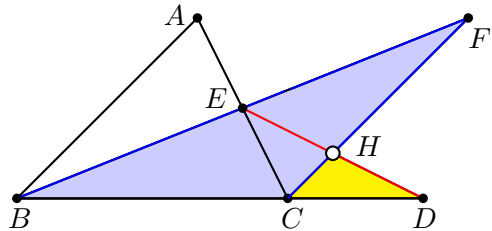

\center{\OneSixteenInnerFigure}
\caption{Proof of I.16 using inner Pasch.}
\label{figure:OneSixteenOne}
\end{figure}

It is instructive to see how inner and outer Pasch are needed to provide
infrastructure for Euclid.  We illustrate with Proposition~I.16, the exterior
angle inequality. That proposition says, 

\begin{quote}  In any triangle, if one of the sides be produced, the exterior
angle is greater than any of the interior and opposite angles. 
\end{quote}
  Refer to  Figure~\ref{figure:OneSixteenOne}.  $ABC$ is the triangle and 
$ACD$ is the exterior angle, which is asserted to be greater than angle $BAC$
and greater than angle $ABC$.  To prove that, 
Euclid constructs $F$ with $EF=EB$, and
 proves triangle $AEB$ is equal to triangle $CEF$,
so in order to prove that angle $ACD > BAC$ it suffices to prove
that $ECF < ACD$.  Euclid justified that with Common Notion~5, the 
whole is greater than the part.  But long before Pasch, one might 
have objected,  how do we know that angle $ECF$ actually {\em is}
a part of $ACD$?  That question needs the same answer that the 
modern definition of angle ordering requires:  the construction 
of the point $H$.  
 In the proof we gave, inner Pasch is used to construct $H$,
 as shown in the figure.

\begin{figure}[ht]
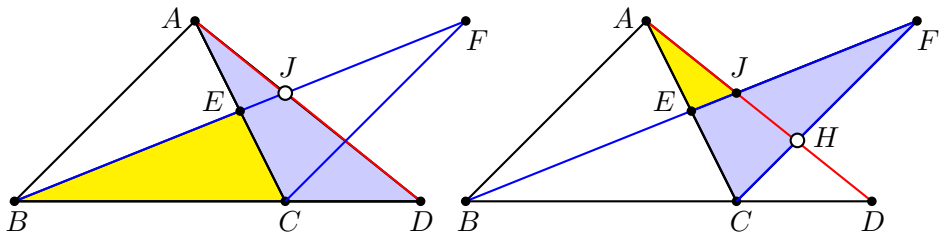

\center{\OneSixteenOuterFigure \OneSixteenFigureTwo}
\caption{Proof of I.16 using outer Pasch.}
\label{figure:OneSixteenTwo}
\end{figure}

One can also prove I.16 from outer Pasch, instead of inner Pasch.
It requires two applications of outer Pasch, as shown in 
Figure~\ref{figure:OneSixteenTwo}.

\section{Angle bisection.}

Euclid~I.9 gives a construction to bisect a given angle $PQR$.  Namely, 
lay off equal segments on the two sides of the angle,  connect their
endpoints $AB$, and use I.1  to construct an equilateral triangle 
$ABC$.  Then $QC$ is the angle bisector.  Ah, but to determine a line,
we must have $Q \neq C$.  And if the original triangle is equilateral,
$C$ {\em will} be $Q$.  So what then?  Well, says the devil, then 
just take the {\em other} intersection point of the two circles in 
Proposition~I.1.  That is, we should modify the circle-circle continuity 
axiom to say there are {\em two} points of intersection of the circle.
All right, let us suppose that is done.   Then we confront the real
difficulty of the proof:  why is $QC$ the angle bisector?  In fact, 
why does it even lie in the plane of angle $PQR$?   

For this proposition to be correct, the definition of the bisector 
of an angle must say that there are points $U$ and $V$ on the 
sides of the angle, and the bisector connects the vertex with some 
point between $U$ and $V$.  Now, even if we assume there are two
equilateral triangles on $AB$, and one of them has a third vertex 
$C$ different from $Q$, there is no apparent reason why $QC$ must meet
$AB$,  as the revised hypothesis requires.

\begin{figure}[ht]
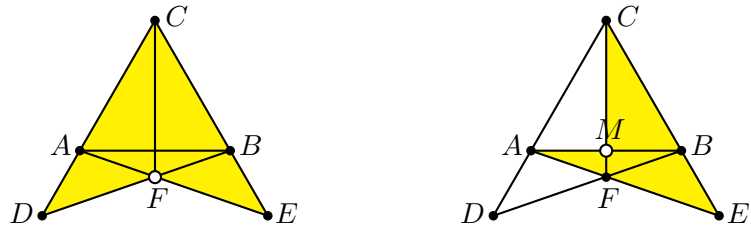

\caption{To bisect angle $ACB$, given  $AC$ equal to $BC$, make
$AE = BD$, and then
construct  $F$  by inner Pasch. One more application of Pasch bisects $AB$.}
\label{figure:proclus}
\center{\ProclusBisectionFigure \ProclusBisectionFigureTwo}
\end{figure}

Proclus already pointed out and attempted to repair 
some of these difficulties in 450 CE,
see  \cite[p.~214]{proclus} and Heath's commentary on I.9 \cite{euclid1956}.
Of course Proclus did not have Pasch's axiom at his disposal; but 
his proof is easily completed using inner Pasch. 
See Figure~\ref{figure:proclus}. (In that figure, we shade
the whole quadrilateral formed by the four points to which inner Pasch
is applied, since the choice of three of them is arbitrary.) 

\FloatBarrier

Proclus observes that triangles $ABE$ and $BAD$ are 
equal (congruent), since angles $EAB$ and $DBA$ are equal by 
Proposition~I.5.  Then $AF$ and $BF$ are equal, by I.6. 
Then triangles $ACF$ and $BCF$ are equal.  Then (by definition of 
equal angles) $CF$ is the desired bisector. 

 After that,
we can proceed
directly from Proclus's proof of I.09 to I.10, since the midpoint
of $AB$ can be constructed with one more application of inner Pasch,
as shown in the second part of Figure~\ref{figure:proclus}.

Proclus also noticed and repaired the problem that it needs to be 
proved that $F$ lies in the interior of the angle.  Since Proclus
did not have inner Pasch available, he made another argument; but 
inner Pasch solves that problem too, as well as the problem of 
showing that the bisector lies in the plane of the angle, 
which neither Proclus nor Heath noticed.%
\footnote{In \cite{beeson2019}, not having studied Proclus enough, we 
used instead a proof from Gupta's 1965 thesis \cite{gupta1965},
which can also be seen as  Lemma~7.25  in \cite{schwabhauser}.  With the 
same figure,  Gupta proves that $M$ is the midpoint of line $AB$.
We used Gupta's proof to prove Euclid~I.10 (line bisection) and 
then used I.10 to prove I.09.  But Gupta's proof is complicated,
because it avoids using circles.  Proclus's proof is simpler, and 
allows us to preserve
Euclid's order of the propositions.}

\section{Two dimensions or three?}
Euclid Books I--IV are commonly thought to be about plane geometry,
but consider: 

\begin{itemize}
\item There is a definition of {\em plane}. 
\item The definition
of {\em parallel} mentions that the two lines must be in the same plane.
\item There is no ``dimension axiom,'' such as Tarski's axiom that three
points each equidistant from points $P,Q$ must be collinear, which 
guarantees that all points lie in a plane.
\item In the last Book, Euclid takes up the Platonic solids, and certainly
uses the results of Book~I. 
\end{itemize}
Nevertheless all the diagrams in Books~I--IV 
appear to be planar figures.  We conclude that Euclid's intention
was to present theorems (and proofs) that are valid in every plane. 
Remember that Euclid did not have our modern conception of ``model'' of 
a geometrical theory.  There was just one true space, and it was three-dimensional,  containing many planes.

But there are several places in Book~I where this seems to have  been 
forgotten.  For example, 
Proposition~7  says that if $ABC$ and $ABD$ are
two triangles with $AC=AD$ and $BC=BD$, and $C$ and $D$ are on the 
same side of $AB$,  then $C=D$.    
The figure Euclid gives is supposed to be impossible,
but as soon as you remember it might be in three dimensions, it looks
very possible.  It is saved from being mistaken by the hypothesis 
mentioning ``same side,''  which forces the diagram to be planar; but 
Euclid did not define ``same side,'' nor did he use that hypothesis 
in the proof.   The definition of ``same side'' is discussed
below; for now we are focussing on the dimension issue.

Euclid's last definition in Book~I is
\begin{quote}
 Parallel straight lines are straight lines which, being in the same plane and being produced indefinitely in both directions, do not meet one another in either direction.
\end{quote}
  
The inclusion of the requirement that the lines be in the same plane
shows that Euclid's omission of a dimension
axiom was not a simple oversight:  he meant to allow 
for the possibility that lines might {\em not} lie in 
the same plane.  But he did 
not define ``lie in the same plane.''   
Whatever Euclid meant by his definition, he found it 
obvious that if $AB$ is parallel to $CD$, then 
$CD$ is parallel to $AB$.  I say that because  this fact is used 
(in the proof of Proposition~I.30)
without even being mentioned (let alone proved). 
This property is called the ``symmetry of parallelism.''

\begin{figure}[ht]
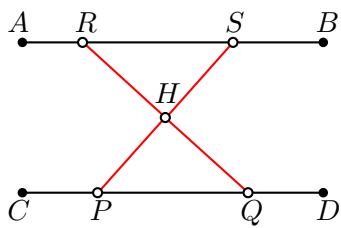

\center{\FigureParallelDefFour}
\caption{$AB$ and $CD$ are coplanar,  as witnessed by  $RSPQH$.}
\label{figure:crisscross}
\end{figure}

We repair this omission by defining two (distinct) lines $AB$
and $CD$ to lie in the same plane if they are linked by a 
``crisscross'' configuration, as shown in Figure~\ref{figure:crisscross}.

With this definition of ``coplanar,''  we can take Euclid's 
definition of ``parallel'' literally:  $AB$ is parallel to 
$CD$ if they are coplanar and do not meet.  This definition makes
the symmetry of parallelism an immediate consequence, and it
also makes it evident that if $AB$ is 
parallel to $CD$, and  any or all of the four points are moved
along their respective lines to new (distinct) positions,  then 
$AB$ is still parallel to $CD$. %
\footnote{
 We define ``Tarski-parallel''
by ``$AB$ and $CD$ do not meet, and $C$ and $D$ lie on the same side of 
$AB$.''   This is clearly not what Euclid intended, as to Euclid it 
seems obvious that if $AB$ is parallel to $CD$ then $CD$ is parallel to 
$AB$, but it requires the plane separation theorem to prove that 
about Tarski-parallel. 
 On the other hand, 
the two definitions can be proved equivalent.  It follows that if $AB$ 
and $CD$ are parallel then $A$ and $B$ are on the same side of $CD$, which 
is quite often actually necessary,  but never remarked by Euclid. 
This is another example of ``infrastructure.'' 
}

\begin{figure}[ht]
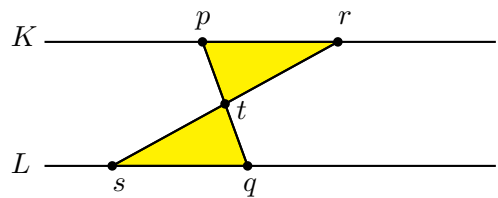
   
\caption{Transversal $pq$ makes alternate interior angles equal with $L$ and $K$,  if $pt=tq$ and 
$rt=st$.}
\label{figure:AlternateInteriorAnglesFigure}
\hskip 2.5cm
\AlternateInteriorAnglesFigure
\end{figure}

Euclid's Postulate~5,  the ``parallel postulate,''  mentions the 
concepts of alternate interior angles, and the concept that the 
two angles on the same side of a transversal ``make together'' more 
than or less than two right angles.  The intention clearly is that 
the lines involved all lie in the same plane,  which will have to 
somehow involve intersection points of some lines.   
It is sometimes possible to reduce theorems about angles directly
to statements about points and the equality relation between 
segments. In particular, it is not 
necessary to develop the theory of angle ordering to state Euclid's parallel postulate.
In Figure~\ref{figure:AlternateInteriorAnglesFigure}, we show how to 
translate the concept ``equal alternating interior angles'' into the
formal language we used.

\smallskip    

Another place where Euclid forgot about three dimensions is 
in Proposition~I.9, the bisection of an angle.  This proof was
discussed above, but here we take up the 3-dimensional aspect
of it.  The problem is that the equilateral triangle constructed
by Proposition~I.1 should introduce a new point by its vertex, not 
just give back the original angle's vertex. 
One may wish to take ``the other'' equilateral triangle, but 
Proposition~I.1 only says there is one, not two.  But in the 
absence of a dimension axiom, there is a whole
circle $C$ of intersection points,  on the plane that bisects $AB$.
In the absence of a dimension axiom,
{\em we have to think of circles as spheres}.  The circle-circle 
axiom is still valid, but even if we assume the two intersection points are on a diameter 
of that circle $C$, it might be tilted out of the plane of $PQR$.
I do not say that Euclid had that mental picture; only that he 
did not have a dimension axiom, and apparently quite purposely, 
so that we who formalize his work must remember that
there is no dimension axiom. 

Euclid does not define 
``rectangle.''  One would like to define it as a quadrilateral with four right 
angles.  It is a theorem that such a figure must lie
in a plane.  However,  the proofs we found involve reasoning ``in three dimensions.''
Even though Euclid Book I has no dimension axiom, and we must therefore 
be careful not to assume one, nevertheless all the {\em proofs} in Book~I 
deal with planar configurations.   We therefore define ``rectangle'' to be a quadrilateral with four 
right angles, whose diagonals cross,  that is, meet in a point.  This condition 
is one way of specifying that a rectangle lies in a plane.  We can then prove 
that a rectangle is a parallelogram.  

Euclid defines a square to be a quadrilateral with at least one
 right angle, in which 
all the sides are equal.
But in I.46 and I.47 the proofs work as if the definition 
required all four angles to be right, so we take that as the definition.
He does not specify that all four vertices lie in 
the same plane.  This is not trivial to prove, but we did prove it, so Euclid's
definition does not require modification.   

\section{Sides of a line and the crossbar theorems.}
We have already discussed Proposition~I.7, which mentions the undefined
``same side'' but never uses it in the proof. 
Since Euclid never defined {\em same side},
there is no obvious way to fix it.   Hilbert worked in plane geometry
in the strong sense, so he did not need to define {\em same side} 
in a way that works in space.  

\begin{figure}[ht]
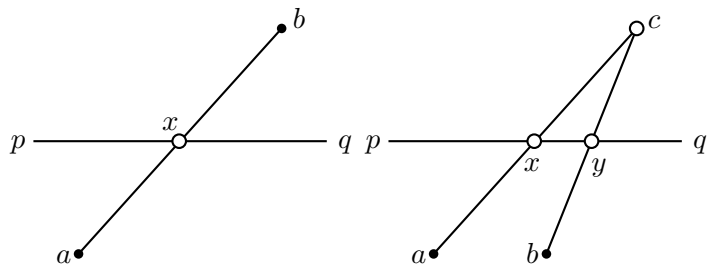

\center{\OppositeSideFigure\quad \SameSideFigure}
\caption{(Left) $a$ and $b$ are on the opposite side of $pq$.
(Right) $a$ and $b$ are on the same side of $pq$ if there exist
points $x$ and $y$ collinear with $pq$, and a point $c$,
 such that $\B(a,x,c)$ and $\B(b,y,c)$.}
\label{figure:sides}
\end{figure}

That notion was, apparently, first defined by M. Pasch in 1882
 \cite[p.~27]{pasch1882}, but only under the assumption 
that the points lie in the same plane.  To remove co-planarity as a 
primitive notion from the definition was first done by Tarski (as far
as I can determine).  He defined two points $a$ and $b$ to be on opposite sides
of $pq$ if there is a point between $a$ and $b$ collinear with $pq$,
and defined $a$ and $b$ to be on the same side of $pq$ if they are
both on the opposite side of $pq$ from the same point $c$.
(See Figure~\ref{figure:sides}.)

Once these concepts are defined, one can use (both inner and outer)
Pasch to prove the {\em plane separation theorem}:  if
$C$ and $D$ are on the same side of $AB$, and 
$D$ and $E$ are on opposite sides of $AB$, then $C$ and $E$
are on opposite sides.  Since neither Pasch nor ``same side'' 
occurs in Euclid, this is not a Euclidean theorem; it is 
infrastructure provided by Tarski and Szmielew, 2300 years later.
Yet we needed it to correct the proofs of not only Euclid~I.7, but 
also Propositions~11, 27, 28, 29, 30, 35, 42, 44, 45, 46, and 47.
These corrections would still be required even if we did add a 
dimension axiom, as,  even in the plane, a line has two sides and 
the plane separation theorem requires a proof.  

\begin{figure}[ht]
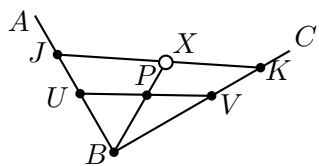

\center{\CrossbarFigure}
\caption{Crossbar theorem: drawing $BP$ through to $JK$.}
\label{figure:crossbar}
\end{figure}

A similar piece of infrastructure, also closely related to 
Pasch's axiom, is the crossbar theorem.  See Figure~\ref{figure:crossbar}.

 This theorem says that 
if $P$ is a point in the interior of angle $ABC$ (which means that 
$P$ lies between two points $U$ and $V$ on the rays $BA$ and $BC$,
respectively),  and if $J$ and $K$ are any two points on those rays,
then the ray $BP$ will meet the ``crossbar'' $JK$ in some point $X$.
If Euclid needs such a point, he simply says that $BP$ is ``drawn through''
to $X$.    One 
place where the crossbar theorem is needed is to prove the 
uniqueness of the angle bisector.
To the modern mathematical eye, having proved existence of angle 
bisectors, the next question should be the uniqueness of the angle 
bisector.  This never occurs as a proposition in Euclid,  and its
omission causes no harm in Book~I.  But it is definitely required 
to fix the proof of III.20,  and it takes more than 110 inferences,
most of which are ``infrastructure'' steps,
concerning the collinearity or noncollinearity of points, the equality 
of various angles,  the transitivity of the less-than relation on angles, 
etc.  These are all things that Euclid usually did not mention.

There are  several versions of the crossbar theorem.  More than
one version is needed, because sometimes we need to know the order
of the points $BPX$.  If we assume $\B BUJ$ and $\B BVK$, then 
 the conclusion can be 
strengthened to $\B BPX$;  another version has those betweenness relations
reversed.  One is proved with two applications of inner Pasch, the other 
with two applications of outer Pasch.  But sometimes we do not know 
the order of the points $BUJ$, so we need the version stated with rays, too.
  
All this infrastructure is required because of the modern insistence on 
requiring the existence of points to be proved,  rather than producing 
the required points by ``drawing through.''

\section{Euclid~4:  all right angles are equal} \label{section:Euclid4}
\begin{figure}[ht]
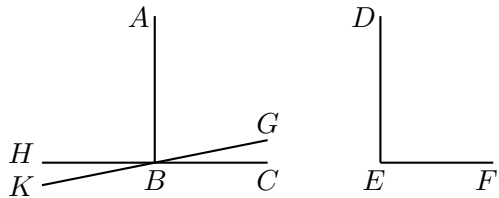

\caption{Proclus's proof of I.4.}
\label{figure:proclusI.4}
\center{\ProclusPostulateFourFigure}
\end{figure}

Over the centuries there were many claims that Euclid~4 is a theorem,
and hence should not be taken as a postulate.  For example, we find 
a proof already in Proclus (450 CE) \cite[p.~148]{proclus}, illustrated
in Figure~\ref{figure:proclusI.4}.

Assume that $ABC$ and $DEF$ are right angles. Proclus says,
``If $DE$ be made to coincide with $AB$, the line $EF$ will fall within
the angle, say, at $BG$.''  Then $H$ and $K$ are constructed,  and
we have $ABC = ABH < ABK = ABG < ABC$,  so $ABC < ABC$, which is 
impossible.  It is the first step that is problematic:  the proof 
appears to depend on the invariance of angles under a rigid motion,
the same flaw that bedevils Euclid's ``proof'' of SAS in I.4.  The 
remedy for that is the angle-copying Proposition I.23.  But 
Euclid appealed to Postulate 4 in proving I.23.  Thus Proclus's 
proof is not correct as it stands.  Moreover, the proof culminates
by saying it is impossible that $ABC < ABC$.  Euclid never proved it;
apparently he considered it to be part of the common notion
``the part is less than the whole.''
 In our formalization,
this theorem is called {\em angletrichotomy2}.  And its (rather lengthy
at 463 lines)
proof uses I.23.%
\footnote{Hilbert smuggled this principle into his axiom system 
without explicit mention, by postulating the uniqueness of the copied
angle in his angle-copying axiom.}
  So there is a {\em second} circularity in Proclus's
proof,  from the modern point of view.  Our formalization 
followed Szmielew's proof, probably discovered at Berkeley in the 1960s,
but published only in 1983 \cite{schwabhauser}; see especially Satz~10.12.
The idea of Szmielew's proof is not so different from Proclus's,
but she replaced the illegal rigid-motion argument with a careful study
of isometries, including reflections in a point or in a line.  Rotations
can be built from reflections in lines.  Of course, this has to be done 
without Postulate~4.  Szmielew's Satz~10.12 is this:  {\em If two right angles have corresponding legs
equal, then the hypotenuses are also equal.}  That theorem 
is easily seen to be equivalent to Postulate~4.  The idea of the proof
is to construct an isometry that takes the corresponding legs onto 
each other.  First a translation brings the two vertices together,
say at point $b$.
Then a rotation makes one leg coincide with the corresponding leg.
This amounts to a correct formalization of Proclus's first step, 
without a circularity.


In \cite{beeson2019}, we followed Szmielew, but in doing so, used Euclid's
construction of perpendiculars based on line-circle and circle-circle continuity.
Szmielew did it without circles, using Gupta's circle-free construction
of perpendiculars from his thesis \cite{gupta1965}.

  \section{Equal figures.} 

The word {\em area} almost never occurs in Euclid's {\em Elements}, despite the 
fact that area is clearly a fundamental notion in geometry. Instead, Euclid speaks of ``equal figures.'' 
 Apparently a ``figure'' is 
a simply connected polygon, or perhaps its interior.  The notion is 
neither defined nor illustrated by a series of examples; for example, it 
is never made clear whether a figure has to be convex,  or even whether 
a circle is a figure, or whether a figure has an interior, or is just
made of lines.  

The notion of ``equal figures'' plays a central role in Euclid.  For 
example, the culmination of Book~I is the Pythagorean theorem.    
Nowadays we would, if required to express the theorem without algebraic
formulas, say that given a right triangle,
 the area of the square on the 
hypotenuse is the sum of the areas of the squares on the sides.  But 
Euclid said instead, that the square on the hypotenuse is equal to the
squares on the sides, taken together.  His proof shows how the two squares
can be cut up into  pieces that can be rearranged to make this equality of figures evident,
given earlier propositions about equal figures. 

Nor was Euclid alone in avoiding the word ``area.''
A century later when Archimedes calculated the area of a circle,
he did not express his result by saying that the area of the circle is 
$\pi$ times the square of the radius.  Instead, he said that circle is 
equal to the rectangle whose sides are the radius and half the circumference.
 (So a circle did count as a figure for Archimedes!) 
 
Why did Euclid avoid the word {area}?  Not because he did not know 
that area can be measured;  it must have been for more abstract, mathematical 
reasons.   Let us consider his problem:
if he were to use the word, he would either have to {\em define} it, 
or put down some {\em postulates} about it.  Both choices offer some 
difficulties.  Area involves assigning a {\em number} to each figure, to 
measure its area.  It is therefore not a purely geometric concept.  
Moreover, even if one is willing to introduce numbers,  that just pushes
the problem back one step:  one must then define or axiomatize numbers.

Euclid's proofs, starting from I.35,  use the notion of 
``equal figures'' without either definition or explicit axiomatization.
He allows himself to paste equal triangles onto equal figures, 
concluding that the results are equal, and justifies it by the common 
notion  {\em if equals be added to equals, the 
wholes are equal}.    He allows himself to cut 
equal triangles off of equal figures, and justifies it by {\em if equals be subtracted from equals, the 
remainders are equal.} If,  with a modern eye, we interpret ``equal figures''
to mean ``figures with equal area,'' these properties look like the 
additivity of area.

  Common Notion 5, ``the whole is greater than 
the part,''  could be taken to imply that a figure cannot be equal to 
a part of itself, and Common Notion 4, ``things which coincide with 
one another are equal to one another,''  could be interpreted to imply
that congruent figures are equal.  
Actually,  Euclid needed one more property: halves of equal
figures are equal, used in Proposition~I.39.  The step that (implicitly) uses that property
occurs in Euclid's text without justification.  

\begin{figure}[ht]
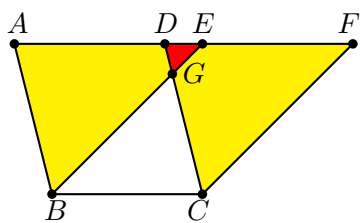

\caption{Euclid's proof of I.35.}
\FigureOneThirtyFiveColored
\label{figure:I.35colored}
\end{figure}

We will give an example of how Euclid reasoned about equal figures,
namely  Euclid I.35.  See Figure~\ref{figure:I.35colored}.

Euclid wants to prove the parallelograms $ABCD$ and $BCFE$
are equal.  He proves the triangles $ABE$ and $DCF$ are 
congruent.  Implicitly, he assumes $DEG$ and $DGE$ are
equal figures (that is, the order of listing the vertices 
does not matter).  Then ``subtracting equals from equals,'' the 
yellow quadrilaterals are equal.  Then, ``adding equals to equals,''
he adds triangle $BCG$ (implicitly assuming $BCG$ is equal to $BGC$)
to arrive at the desired conclusion.

Later generations of mathematicians were not willing to accept Euclid's
over-liberal interpretation of the common notions in support of 
``equal figures.''  See the summary 
discussion with many references \cite[pp.~327--328]{euclid1956}.  In
particular, once mathematicians had some experience with axiomatization,
it became obvious that ``equal figures'' is not a special case of 
equality, since equal figures cannot be substituted for each other 
in every property.   Instead, it is a new relation, and the original
choice that Euclid finessed faced us directly when we wanted to 
proof-check Euclid:  we had either to define
or to axiomatize the notion.   We chose to axiomatize it.%
\footnote{We do not take space here to describe attempts to define it 
by Hilbert and others.  See \cite{beeson2020,hartshorne}
 for a full discussion; \cite{beeson2020} also gives a definition 
 that ``Euclid could have given.''}
Following 
the lead of Hartshorne \cite{hartshorne}, we wrote down fifteen axioms for the 
three primitive notions of ``equal triangles'' and ``equal quadrilaterals''
and ``triangle equal to quadrilateral.''  No figures with more than four 
sides occur in Book~I, so that was sufficient.  These were ``cut-and-paste''
axioms as described above, plus two axioms (first invented
by de Zolt, see~\cite{hartshorne}) saying that if you cut off 
something, the result is not equal to what you had before the cut.  

\section{Making Formal Proofs Readable.}

A proof has two purposes:  to establish beyond doubt that a theorem 
is correct, and to communicate to a human reader {\em why} it is correct.
Our formal proofs improve on Euclid in the first respect, but in 
the second respect they need improvement.  We address that issue now. 

We begin by saying something about what it 
is like to formalize Euclid's proofs.  What one finds is that one
needs a large number (more than sixty percent) of ``invisible infrastructure''
steps, proving statements that Euclid would take for granted because
they appear so in the diagram.  Let me give one example.  Often
to verify the hypotheses of some proposition or lemma we wish to apply,
we need to know that 
an angle $ABC$ is equal to itself.  Before we can write that down, 
we have to verify that the three points $A$, $B$, and $C$ are noncollinear,
so that they really do form an angle.  Euclid always takes that for 
granted if it appears so in the diagram,  but it often requires many 
formal steps, using for example the lemma that if $A$, $B$, $C$ are noncollinear
and $U$ and $V$ are distinct and both are collinear with $AB$, then
$B$, $U$, $V$ are noncollinear.  A chain of three or four applications of that 
lemma may be required to verify that a given triangle really is 
a triangle.  After a while, one can do this ``diagram chasing'' 
as fast as one can type,
but it means that the formal proofs are at least double the length of 
Euclid's and contain many uninteresting steps.  Euclid generally takes
for granted statements of collinearity and betweenness that appear obvious
in the diagram. 

To address the problem of too many and too-detailed ``infrastructure'' steps, 
we created a list of ``trivial'' inferences, or more accurately trivial
justifications,  which are to be suppressed
on output.  For example,  the application of lemmas about collinearity
and noncollinearity.   By putting more or fewer justifications on the
``trivial list,''  we can vary the ``step-size'' or ``grain'' of the 
output proofs.  We call the result a ``proof skeleton''.   We found 
that with a suitable list of trivial steps, the proof skeletons contained
Euclid's steps, and only really essential additions (such as uses of 
Pasch's axiom). For example,  in Prop.~I.16, the proof skeleton has 15 inferences,
while the full formal proof has about 120 inferences, most of which are about
collinearity, non-collinearity, and distinctness of points that appear 
distinct in the diagram.      

A second obstacle to human readability is the fact that 
formal proofs contain only symbols (no words). 
Euclid used no symbols except names for points; so 
if we want to compare our proofs to Euclid's,  we ought to write them 
in Euclid's way.   We 
therefore experimented with machine-generating English-language proofs%
\footnote{Greek-language proofs might have been more authentic.  Output
in any human language can be generated once the ``stock phrases'' are translated.}
from our formal proofs,  ``in the style of Euclid.'' 
 Euclid's natural language consists of a small number
of stock phrases that connect the assertions.  We wrote a script 
that translates proof skeletons (or proofs)  into \LaTeX\ code for an English
version of the proof.  The script follows Euclid by 
prefacing each construction of a new point with a sentence about how it
is constructed. For example, it inserts a 
phrase ``Let $AC$ be bisected at $E$'',  when Prop.~I.10 (the midpoint theorem)
is about to be applied.    
That script first 
extracts a proof skeleton  by deleting trivial steps, and then  
produces a \LaTeX\ document, translating the proof into
``Euclidean English.'' 
In a few seconds, it  processes all 245 formal proofs,
producing  a computer-generated version of Euclid Book~I.
Here we present just two samples from this document.  The reader
is invited to compare them with Euclid's proofs.  In particular, the 
proof of I.16 supplies the  application of Pasch's axiom that is missing
from Euclid's proof; and the proof of I.20 illustrates the formal treatment
of angle ordering, appealing to a definition rather than Common Notion~5.
These proofs demonstrate, I 
assert,  that the translation from formal proofs to human-readable proofs,
while it may be formidable for humans,  is trivial for computers.%
\footnote{The other direction, from human-readable proofs to formal proofs,
is far from trivial.}  
The 
programs that make these transformations are short and simple, and were 
easy to write.

\section*{  Proposition 20 }
{\em In any triangle two sides taken together  
are greater than the third side.}
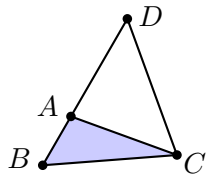
\begin{figure}[ht]
\center{\OneTwentyFigure}
\end{figure}
\medskip

Let $ABC$ be a triangle. 

 It is required to show that $BA$, $AC$ are together greater than $BC$. 
\medskip

\noindent
 \hspace{0pt}Let $BA$ be produced in a straight line to $A$, making $AD$ equal to $CA$.  Then

\hspace{0pt}$A$ is between $B$ and $D$ and $AD$ is equal to $CA$.\justification{extension}
\hspace{0pt}Triangle $ADC$ is isoceles with base $DC$.\justification{defn:isosceles}
\hspace{0pt}Angle $ADC$ is equal to angle $ACD$.\justification{I.5}
\hspace{0pt}Angle $ADC$ is equal to angle $DCA$.\justification{equalanglestransitive}
\hspace{0pt}Angle $DCB$ is greater than  angle $ADC$.\justification{defn:anglelessthan}
\hspace{0pt}Angle $ADC$ is equal to angle $BDC$.\justification{defn:equalangles}
\hspace{0pt}Angle $DCB$ is greater than  angle $BDC$.\justification{angleorderrespectscongruence2}
\hspace{0pt}Angle $DCB$ is greater than  angle $CDB$.\justification{angleorderrespectscongruence2}
\hspace{0pt}Angle $BCD$ is greater than  angle $CDB$.\justification{angleorderrespectscongruence}
\hspace{0pt}$BD$ is greater than $BC$.\justification{I.19}
\qed

\section*{Proposition 16}
{\em In any triangle, if one of the sides be produced,
the exterior angle is greater than either of the interior
and opposite angles.}

\begin{figure}[ht]
\center{\OneSixteenInnerFigure}
\label{figure:OneSixteenRepeat}
\end{figure}

Let $ABC$ be a triangle, and let one side of it $BC$ be produced
to $D$; then $C$ is between $B$ and $D$. \\

It is required to show that the exterior angle $ACD$ is greater than 
the interior and opposite angle  $BAC$.%
\footnote{\,I.16 asserts that the exterior angle is greater than both
interior angles.  Like Euclid, we here present only the proof for one exterior angle.}
\medskip

\noindent
\hspace{0pt}Let $AC$ be bisected at $E$.  Then

\hspace{0pt}$E$ is between $A$ and $C$ and $EA$ is equal to $EC$.\justification{I.10}
\hspace{0pt}Let $BE$ be produced in a straight line to $E$, making $EF$ equal to $EB$.  Then

\hspace{0pt}$E$ is between $B$ and $F$ and $EF$ is equal to $EB$.\justification{extension}
\hspace{0pt}Let $AC$ be produced in a straight line to $C$, making $CG$ equal to $EC$.  Then

\hspace{0pt}$C$ is between $A$ and $G$ and $CG$ is equal to $EC$.\justification{extension}
\hspace{0pt}Angle $BEA$ is equal to angle $CEF$.\justification{I.15}
\hspace{0pt}Angle $AEB$ is equal to angle $CEF$.\justification{equalanglestransitive}
\hspace{0pt}$AB$ is equal to $CF$, and Angle $EAB$ is equal to angle $ECF$, and Angle $EBA$ is equal to angle $EFC$.\justification{I.4}
\hspace{0pt}Angle $BAC$ is equal to angle $BAE$.\justification{equalangleshelper}
\hspace{0pt}Angle $BAC$ is equal to angle $EAB$.\justification{equalanglestransitive}
\hspace{0pt}Angle $BAC$ is equal to angle $ECF$.\justification{equalanglestransitive}
\hspace{0pt}Angle $ECF$ is equal to angle $ACF$.\justification{equalangleshelper}
\hspace{0pt}Angle $BAC$ is equal to angle $ACF$.\justification{equalanglestransitive}
\hspace{0pt}Let $CF$ and  $ED$ meet at $H$.  Then

\hspace{0pt}$H$ is between $D$ and $E$ and $H$ is between $F$ and $C$.\justification{Pasch-inner}
\hspace{0pt}Angle $BAC$ is equal to angle $ACH$.\justification{equalangleshelper}
\hspace{0pt}Angle $BAC$ is equal to angle $ACF$.\justification{equalangleshelper}
\hspace{0pt}Angle $BAC$ is equal to angle $ACH$.\justification{equalanglestransitive}
\hspace{0pt}Angle $ACD$ is greater than  angle $BAC$.\justification{defn:anglelessthan}
\qed

In these examples, only the italicized informal statement at the top 
and the diagram are human-generated.  The rest, including all the English, all 
the references, and the typesetting, 
is machine-generated.  At last, we have achieved,
and certified by computer, 
the goal that Gerolamo Saccheri stated in the title of his  
1733 book, {\em Euclid Vindicated from Every Blemish}.%
\footnote{Other translators have chosen {\em Euclid Freed of 
Every Flaw}. The title above is the one chosen by de Risi. \cite{saccheri-deRisi}
}

\section{Euclid Vindicated.}
Table~\ref{table:1} compares Euclid's {\em Elements}
(Books I to IV) with the changes we made in formalizing Euclid. 

\begin{table}[ht]
\caption{Changes to Propositions and Axioms.}
\label{table:1}
\begin{center}
\begin{tabular}{l|l|l}
Issue & Euclid &  Changes made \\
\hline
Postulate IV  &long thought provable  &  proved by Szmielew  \\
Proposition I.4  & rejected since antiquity & replaced by 5-line axiom \\
definition of ``parallel''      & ``in the same plane'' & supplied definition   \\
Postulate V    &  ``alternate interior angles'' & supplied definition \\
 connectivity axiom  & missing               &  added \\
betweenness notion     &  missing &            added \\
betweenness axioms   &   missing & identity, symmetry,\\
                     &                       & and transitivity \\
betweenness, basic theorems &  missing  & proved  \\
definition of ``same side''  &  missing  & added Tarski's definition \\
2 or 3 dimensions? & valid in any plane & use Tarski's ``same side'' \\
Pasch & missing axiom &  inner and outer Pasch \\
line-circle and circle-circle & missing axioms             &  added  both \\
common notions for lines   &         & 3,5 proved; 4 dropped\\
equality of angles        & primitive notion &   defined notion \\
less than for angles      & primitive notion (?) &   defined notion \\
common notions for angles &                   & proved, instead \\
                       &                      &of assumed \\
equal figures     & no definition or axioms  &  added 15 axioms \\
rectangle definition & omitted   & four right angles \\
                     &           & and diagonals meet 
\end{tabular}
\end{center}
\end{table}%
Table~\ref{table:2} shows some of the propositions in Book~I 
that needed corrections.
We do not include as ``corrections'' the provision of ``infrastructure''
steps about collinearity and noncollinearity, nor elementary reasoning
about betweenness, nor the many cases where Pasch was needed,
or some lemmas had to be proved, nor proofs where Euclid 
treated only one of several cases (for example I.35).

\begin{table}[ht]
\caption{Corrections to Proofs (refer to Euclid's diagrams).}
\label{table:2}
\begin{center}
\begin{tabular}{p{13pt} l|l|l}
Prop.& Description & Difficulty &  Correction \\
\hline 
I.1& equilateral triangle &  existence of $C$  & circle-circle \\
I.1& equilateral triangle        & $ABC$ might be collinear         & connectivity axiom \\
I.4& SAS   &  superposition   &  5-line axiom \\  
I.7& triangle uniqueness & same side not defined & use Tarski's definition \\
I.7&  triangle uniqueness                    & angle trichotomy  & proved as theorem \\
I.9& angle bisection  & $A$ and $F$ might coincide & use Proclus's proof \\
I.12&  dropped perpendicular & Why do $G$ and $E$ exist? & line-circle \\
I.16&  exterior angle        & Why is $ECD > ECF$?   & Pasch \\
I.22&  triangle construction & why does $K$ exists?  & circle-circle \\
I.22&  triangle construction & why does $DE$ meet circles?& line-circle \\
I.27&  parallel construction  & ``alternate angles'' undefined & $AD$ 
and $EF$ \\
  & & &  must meet \\
I.32& exterior angle & see I.16 & \\
I.33& parallelogram constr. & ``same direction'' & diagonals must meet\\
I.35& parallelograms       & Why is $DEG = EGD$?   & an equal-figures axiom \\
I.35& parallelograms                            & several other steps  & equal-figures axioms\\ 
I.39& equal triangles         & same-side  undefined, unused & 
                                  use Tarki's definition \\
I.46& square definition & definition doesn't match use  
                              & changed definition 
\end{tabular}
\end{center}
\end{table} 
\FloatBarrier

\section{Conclusions.}
The two most characteristic features of Euclid are   geometrical 
diagrams, and chains of 
logical reasoning about those diagrams.  The exact relation between these two features
has been a concern of every thoughtful reader of Euclid, right from the beginning.  
The reasoning is guided by the diagram;  but sometimes it is led astray by the diagram,
too!  It took two millennia   to {\em separate} the two features,
inventing {\em symbolic logic}--meaningless chains of symbols representing correct
inferences made according to precise rules.  As the logician J. Barkley Rosser expressed
it \cite[p.~7]{rosser1978}:

\begin{quote} This does not mean that it is now any easier to discover a proof for a difficult
theorem.  This still requires the same high order of mathematical talent as before.
However, once the proof is discovered, and stated in symbolic logic, it can be checked
by a moron.
\end{quote}

It can even be checked by a computer. Rosser's teacher, Kleene, was once asked why he 
wrote his proofs so formally.  He replied, ``How else can I be sure they are right?'' 
And that is the most obvious, and most important,
 result of proof-checking Euclid:  Now we are {\em sure} that
the proofs are correct.  Though some of the gaps and errors we uncovered were known 
for a long time, others were not, so mere human checking did not really do the job.

One striking feature of these formal proofs is that they are longer than the proofs 
mathematicians write,  usually by a factor of about four.
I have called those extra steps ``infrastructure.''  In Euclid they are mostly about 
collinearity, noncollinearity, and betweenness.  They represent facts that a human reader infers from
the diagram and takes for granted without explicit proof.  Even 
when we {\em think}  we are checking a proof carefully, we are skipping many necessary small
steps, jumping over those steps to reach a conclusion that we believe {\em on some other grounds}, for example, on the appearance of a diagram.
When we see that the evidence of our eyes or intuition (the diagram) is confirmed by the successful 
completion of a long chain of logical reasoning, that produces a feeling of 
satisfaction that is the heart of mathematics.  That's why eleven-year old Bertrand 
Russell called Euclid ``delicious.''     Neither a diagram without reasoning,
nor a meaningless chain of inferences, deserves the name, ``mathematics.''   

\begin{acknowledgment}
{Acknowledgments}
I am indebted to 
John Baldwin,
Pierre Boutry,
Erwin Engeler,
Vincenzo de Risi,
Julien Narboux,
Victor Pambuccian,
Dana Scott, 
Albert Visser,
and 
Freek Wiedijk,
for many conversations and emails about Euclid and formalization.
I dedicate this paper {\em in memoriam} to Marvin Jay Greenberg,
who first introduced me to axiomatic geometry.  
\end{acknowledgment}


\begin{thebibliography}{10}
\expandafter\ifx\csname urlstyle\endcsname\relax
 \providecommand{\doi}[1]{doi:\discretionary{}{}{}#1}\else
 \providecommand{\doi}{doi:\discretionary{}{}{}\begingroup
  \urlstyle{rm}\Url}\fi

\bibitem{beeson2020}
Beeson, M. (2020).
\newblock On the notion of equal figures in {E}uclid.
\newblock Available at: ArXiv 2008.12643, math.LO

\bibitem{beeson2019}
Beeson, M., Narboux, J., Wiedijk, F. (2019).
\newblock Proof-checking {E}uclid.
\newblock \textit{Ann. Math. Artif. Intell.}  85(2):
  213--257.
\newblock \doi{10.1007/s10472-018-9606-x}.

\bibitem{beeson2017a}
Beeson, M., Wos, L. (2017).
\newblock Finding proofs in {T}arskian geometry.
\newblock \textit{J.~Automat. Reason.}, 58(1): 181--207.

\bibitem{boutryphd}
Boutry, P. (2018).
\newblock {On the Formalization of Foundations of Geometry}.
\newblock Ph.D. thesis. University of Strasbourg.

\bibitem{narboux2017b}
Boutry, P., Braun, G., Narboux, J. (2017).
\newblock Formalization of the arithmetization of {E}uclidean plane geometry
  and applications.
\newblock \textit{J.~Symbolic Comput.}, 90(1): 149--168.

\bibitem{narboux2017c}
Boutry, P., Gries, C., Narboux, J., Schreck, P. (2019).
\newblock {Parallel postulates and continuity axioms: a mechanized study in
  intuitionistic logic using Coq}.
\newblock \textit{{J.~Automat. Reason.}}, 62(1): 1--68.\\
\newblock \doi{10.1007/s10817-017-9422-8} 

\bibitem{narboux2017}
Braun, G., Narboux, J. (2017).
\newblock A synthetic proof of {P}appus' theorem in {T}arski's geometry.
\newblock \textit{J.~Automat. Reason.}, 58(2): 209--230.
\newblock \doi{10.1007/s10817-016-9374-4} 

\bibitem{deRisi2016b}
De~Risi, V. (2016).
\newblock The development of {E}uclidean axiomatics.
\newblock \textit{Archive for History of Exact Sciences}. 70(6): 591--676.
\newblock \doi{10.1007/s00407-015-0173-9} 

\bibitem{deRisi2019}
De~Risi, V. (2019).
\newblock Leibniz on the continuity of space.
\newblock In: V.~De~Risi, ed., \textit{Leibniz and the Structure of Sciences:
  Modern Perspectives on the History of Logic, Mathematics, Epistemology}.
  Cham: Springer International Publishing, pp. 111--169.
\newblock \doi{10.1007/978-3-030-25572-5_4} 

\bibitem{deRisi2020}
De~Risi, V.  (forthcoming).
\newblock {E}uclid\textquoteright{}s common notions and the theory of
  equivalence.
\newblock \textit{Foundations of Science}: 1--24.
\newblock \doi{10.1007/s10699-020-09694-w} 

\bibitem{charlesriver}
Editors, C.~R. (2014).
\newblock \textit{The Library of {A}lexandria: The History and Legacy of the
  Ancient World's Most Famous Library}.
\newblock Boston: CreateSpace Independent Publishing Platform.

\bibitem{el-abbadi2020}
El-Abbadi, M. (2020).
\newblock Library of {A}lexandria.
\newblock \textit{Encyclopedia Brittanica}.

\bibitem{simson}
Euclid (1787).
\newblock \textit{The {E}lements of {E}uclid, viz. the first six books, together
  with the eleventh and the twelfth}.
\newblock Edinburgh: Nourse and Balfous, 7th ed.
\newblock (Simson, R., trans.)  Available from Bibliotheque Nationale at
 gallica.bnf.fr/ark:/12148/bpt6k1163221v

\bibitem{euclid1956}
Euclid (1956).
\newblock \textit{The Thirteen Books of The {E}lements}.
\newblock (Heath, T.~L., trans.)
\newblock New York: Dover.
\newblock Three volumes. 

\bibitem{gelernter1959}
Gelernter, H. (1963).
\newblock Realization of a geometry theorem-proving machine.
\newblock In:  Feigenbaum, E., J.~Feldman, J., eds., \textit{Computers and Thought}.
  New York: McGraw-Hill, pp. 134--152.

\bibitem{gelernter1960}
Gelernter, H., Hansen, J.~R., Loveland, D.~W. (1963).
\newblock Empirical explorations of a geometry-theorem proving machine.
\newblock In: Feigenbaum, E., J.~Feldman, J., eds., \textit{Computers and Thought}.
  New York: McGraw-Hill, pp. 153--167.

\bibitem{gupta1965}
Gupta, H.~N. (1965).
\newblock \textit{Contributions to the Axiomatic Foundations of Geometry}.
\newblock Ph.D. thesis. University of California, Berkeley.
Available through ProQuest, LLC.

\bibitem{hartshorne}
Hartshorne, R. (2000).
\newblock \textit{Geometry: Euclid and Beyond}.
\newblock New York: Springer.

\bibitem{heath-vol1}
Heath, S.~T. (1921).
\newblock \textit{A History of {G}reek Mathematics, vol. I: From {T}hales to
  {E}uclid.}
\newblock Oxford: Clarendon Press.

\bibitem{hilbert1899}
Hilbert, D. (1960).
\newblock \textit{Foundations of Geometry ({G}rundlagen der {G}eometrie)}.
\newblock La Salle, IL: Open Court.
\newblock Second English edition, translated from the tenth German edition by
  Leo Unger. Original publication date, 1899.

\bibitem{hirsch}
Hirsch, D., van Haften, D. (2010).
\newblock \textit{{A}braham {L}incoln and the Structure of Reason}.
\newblock New York: Savas Beatie LLC.

\bibitem{kennedy2002}
Kennedy, H. (2002).
\newblock \textit{Life and Works of {G}uiseppe {P}eano}.
\newblock San Francisco: Peremptory Publications.

\bibitem{ketcham}
Ketcham, H. (1901).
\newblock \textit{The Life of Abraham Lincoln}.
\newblock New York: The Perkins Book Company.

\bibitem{mollerup1904}
Mollerup, J. (1904).
\newblock Die {B}eweise der ebenen {G}eometrie ohne {B}enutzung der
  {G}leichheit und {U}ngleichheit der {W}inkel.
\newblock \textit{Math. Ann.} 58: 479--496.

\bibitem{oleary}
O'Leary, D.~L. (2001).
\newblock \textit{How Greek Science Passed to the {A}rabs}.
\newblock New Delhi: Goodword Books.

\bibitem{pasch1882}
Pasch, M. (1882).
\newblock \textit{{V}orlesung \"uber {N}euere {G}eometrie}.
\newblock Leipzig: Teubner.

\bibitem{peano1889}
Peano, G. (1889).
\newblock \textit{Principii de Geometria}.
\newblock Torino: Fratelli Bocca.

\bibitem{proclus}
Proclus (1970).
\newblock \textit{A Commentary on the First book of {E}uclid}.
\newblock Princeton: Princeton Univ. Press.

\bibitem{rigby1970}
Rigby, J.~F. (1970).
\newblock Axioms for absolute geometry, {III}.
\newblock \textit{Canadian J.~Math.} 22(1): 185--190.

\bibitem{rosser1978}
Rosser, J.~B. (1978).
\newblock \textit{Logic for Mathematicians}, 2nd ed.
\newblock Mineola, NY: Dover Publications.

\bibitem{russell1902}
Russell, B. (1902).
\newblock The teaching of {E}uclid.
\newblock \textit{Math. Gaz.} 2(33): 165--167.

\bibitem{russell-autobiography}
Russell, B. (2000).
\newblock \textit{Autobiography of Bertrand Russell}.
\newblock London:Routledge.

\bibitem{saccheri-deRisi}
Saccheri, G. (2014).
\newblock \textit{Euclid Vindicated from Every Blemish}.
\newblock Cham: Birkh\"auser (Springer).
\newblock (Halstead, G. B., Allegri, L., trans.; De Lisi, V., ed.); original date 1733.

\bibitem{sarton2}
Sarton, G. (1993).
\newblock \textit{Hellenistic Science and Culture in the Last Three Centuries
  {B}.{C}.}
\newblock Mineola, NY:Dover.

\bibitem{schwabhauser}
Schwabh\"auser, W., Szmielew, W., Tarski, A. (1983).
\newblock \textit{{M}etamathematische {M}ethoden in der {G}eometrie: {T}eil {I}:
  {E}in axiomatischer {A}ufbau der {E}uklidischen {G}eometrie. {T}eil {II}:
  {M}etamathematische {B}etrachtungen ({H}ochschultext)}.
\newblock Berlin:Springer--Verlag.
\newblock Reprinted 2011 by Ishi Press, with a new foreword by Michael Beeson.

\bibitem{strabo}
Strabo (2016).
\newblock \textit{Geography}.
\newblock Hastings, East Sussex, UK: Delphi Classics.
\newblock (Hamilton, H.~C., Falconer, W., trans.)

\bibitem{tarski1959}
Tarski, A. (1959).
\newblock What is elementary geometry?
\newblock In:  Henkin, L.,  Suppes, P.,  Tarksi, A., eds., \textit{The Axiomatic
 Method, with Special Reference to Geometry and Physics. Proceedings of an
  International Symposium held at the {U}niversity of {C}alifornia, {B}erkeley, {D}ecember 
  26, 1957--{J}anuary 4, 1958}. Studies in Logic and the Foundations of
  Mathematics. Amsterdam: North-Holland, pp. 16--29.
\newblock Available as a 2007 reprint, Brouwer Press, ISBN 1-443-72812-8.

\bibitem{tarski-givant}
Tarski, A., Givant, S. (1999).
\newblock Tarski's system of geometry.
\newblock \textit{Bull. Symb. Log.}  5(2): 175--214.

\bibitem{veronese1891}
Veronese, G. (1891).
\newblock \textit{Fondamenti di Geometria a Pi\`u Dimensioni e a Pi\`u Specie di
  Unit\`a Rettilinee Esposti in Forma elementare. Lezioni per la Scuola di
  Magistero in Matematica}.
\newblock Padova: Tipografia del Seminario.

\end{thebibliography}

\vskip 1cm

\noindent
 
 Michael Beeson retired from San Jos\'e State University in 2013.  He 
 has been studying the foundations of geometry since 2006.  He is the 
 author of the software {\em MathXpert}, various papers on automated
 deduction, the foundations of constructive mathematics, minimal surfaces,
 and triangles.  See his web pages  {\url{www.michaelbeeson.com}} and 
 \url{www.helpwithmath.com}.  He can be contacted at {\tt beesonpublic@gmail.com}.

\end{document}